\documentstyle[12pt,amsfonts,amssymb]{amsart}

\newcommand{\beq}{\begin{equation}}
\newcommand{\eeq}{\end{equation}}

\newcommand{\rarrow}{\rightarrow}

\newcommand{\boldZ}{{\Bbb Z}}

\newtheorem{theorem}{Theorem}[section]
\newtheorem{prop}[theorem]{Proposition}
\newtheorem{importnota}[theorem]{Important Notation}
\newtheorem{prblm}[theorem]{Problem}
\newtheorem{notation}[theorem]{Notation}
\newtheorem{defin}[theorem]{Definition}
\newtheorem{caution}[theorem]{Caution}
\newtheorem{remark}[theorem]{Remark}
\newtheorem{lemma}[theorem]{Lemma}
\newtheorem{construction}[theorem]{Construction}
\newtheorem{corollary}[theorem]{Corollary}
\newtheorem{example}[theorem]{Example}
\newtheorem{conclusion}[theorem]{Conclusion}
\newtheorem{triviality}[theorem]{Triviality}
\newtheorem{proto}[theorem]{Prototype Quasifibration}
\newtheorem{cauex}[theorem]{Cautionary Example}
\newtheorem{subth}{ }[theorem]
\newtheorem{ssubth}{ }[subth]

\newcommand\nc{\newcommand}
\nc\tri[1]{\begin{triviality}
\label{#1}}
\nc\prt[1]{\begin{proto}
\label{#1}}
\nc\lem[1]{\begin{lemma}
\label{#1}}
\nc\pro[1]{\begin{prop}
\label{#1}}
\nc\thm[1]{\begin{theorem}
\label{#1}}
\nc\cor[1]{\begin{corollary}
\label{#1}}
\nc\dfn[1]{\begin{defin}
\label{#1}}
\nc\sthm[1]{\begin{subth}
\label{#1}}
\nc\exm[1]{\begin{example}
\label{#1}
\begin{em}}
\nc\plm[1]{\begin{prblm}
\label{#1}
\begin{em}}
\nc\rmk[1]{\begin{remark}
\label{#1}
\begin{em}}
\nc\ntn[1]{\begin{notation}
\label{#1}
\begin{em}}
\nc\cau[1]{\begin{caution}
\label{#1}
\begin{em}}
\nc\imn[1]{\begin{importnota}
\label{#1}
\begin{em}}
\nc\cax[1]{\begin{cauex}
\label{#1}
\begin{em}}
\nc\con[1]{\begin{construction}
\label{#1}
\begin{em}}
\nc\ssthm[1]{\begin{ssubth}
\label{#1}
\begin{em}}
\nc\cnc[1]{\begin{conclusion}
\label{#1}
\begin{em}}

\nc\elem{\end{lemma}}
\nc\eprt{\end{proto}}
\nc\ethm{\end{theorem}}
\nc\ecor{\end{corollary}}
\nc\edfn{\end{defin}}
\nc\esthm{\end{subth}}
\nc\epro{\end{prop}}
\nc\etri{\end{triviality}}
\nc\eexm{\end{em}
\end{example}}
\nc\ermk{\end{em}
\end{remark}}
\nc\eplm{\end{em}
\end{prblm}}
\nc\ecau{\end{em}
\end{caution}}
\nc\ecax{\end{em}
\end{cauex}}
\nc\eimn{\end{em}
\end{importnota}}
\nc\entn{\end{em}
\end{notation}}
\nc\econ{\end{em}
\end{construction}}
\nc\ecnc{\end{em}
\end{conclusion}}
\nc\essthm{\end{em}
\end{ssubth}}
\nc\nin{\noindent}
\nc\la{\longrightarrow}
\nc\qqed{\hfill{$\Box$}}
\nc\ct{{\cal T}}
\begin{document}

\title[Grothendieck duality via homotopy theory]{The Grothendieck
Duality Theorem via Bousfield's Techniques and Brown
Representability}
\author{Amnon Neeman\footnotemark[1]}\addtocounter{footnote}{1}%
\footnotetext{The research was partly
supported by NSF grant DMS--9204940}

\address{Department of Mathematics\\
         University of Virginia\\
         Charlottesville, VA 22903}

\maketitle

\abstract
Grothendieck proved that if $f:X\longrightarrow Y$ is a proper
morphism of nice schemes, then $Rf_*$ has a right adjoint, which
is given as tensor product with the relative canonical
bundle.
The original proof was by patching local data. Deligne proved the
existence of the adjoint by a global argument, and Verdier showed
that this global adjoint may be computed locally.

In this article we show that the existence of the adjoint is
an immediate consequence of Brown's representability theorem.
It follows almost as immediately, by ``smashing'' arguments,
that the adjoint is given by tensor product with a dualising
complex. Verdier's base change theorem is an immediate
consequence.
\endabstract
\tableofcontents

\setcounter{section}{-1}

\section{Introduction}
\label{S0}

Let $f: X \rarrow Y$ be a proper morphism of schemes.  Then, under mild
hypotheses on $f, \ X$ and $Y$, Grothendieck proved that there is a natural
isomorphism
$$
R f_* {\cal RH}om^{}_X(x,f^! y) \simeq {\cal RH}om^{}_Y
(Rf_*x,y)
$$
of objects in the derived category.  We should perhaps briefly remind the
reader
what this means.

Let $D^+(qc/X)$ be the derived category of
bounded--below chain complexes of quasi-coherent sheaves
on $X$.  Let $Rf_*: D^+(qc/X) \rarrow D^+(qc/Y)$ be the right derived
functor of $f_*$.  Then the statement above asserts that $Rf_*$ has
a right adjoint, denoted $f^!$, and that $f^!$ behaves well with respect
to pullbacks by open immersions.

Grothendieck's original proof of the existence of $f^!$ was constructive. It
amounted
to a local computation. Since derived categories are basically unsuited for
local computations, the argument turns out to be quite
unpleasant; see \cite{H}.

There is an abstract way to prove the existence of $f^!$ in the literature.
It is due to Deligne~\cite{D}. The approach was developed and elaborated by
Verdier~\cite{V}, where it is shown that almost everything in \cite{H} can be
obtained
directly from Deligne's result.

All these results assumed the scheme $X$ Noetherian.  Lipman
recently developed a deep theory for
removing the Noetherian hypotheses. The reader is referred
to \cite{L1} and \cite{L2}.

Unfortunately, none of the approaches generalizes
well to $D$--modules.  Given a
morphism $f: X \rarrow Y$, one can define a morphism $Rf_+$ on the
corresponding derived categories of complexes of $D$--modules.
It turns out that $R f_+$ has a right
adjoint.  But the existence of the right adjoint has until now always been
proved by factoring $f$ suitably, defining certain trace maps which
$\grave{a}$ {\em priori} depend on the factorization, and finally proving that
they are independent of factorizations.

What we propose to do here is show that all of the results are direct
consequences of a very simple, general statement about triangulated functors
on triangulated categories.
\bigskip

\nin
{\bf Definition~\ref{compactly generated}} {\em Let ${\cal S}$ be a
triangulated category.  Suppose that all small coproducts exist in
${\cal S}$.  Suppose there exists a set $S$ of objects of
${\cal S}$ such that
\begin{itemize}
\item[(1)] For every $s \in S$, Hom$(s,-)$ commutes with coproducts.
\item[(2)] If $y$ is an object of ${\cal S}$ and Hom$(s,y) = 0$ for all
$s \in S$, then $y = 0$.
\end{itemize}

\nin
Such a triangulated category ${\cal S}$ is called {\em compactly
generated.}}\bigskip

\nin
The Brown Representability Theorem, in one version, states:
\bigskip

\noindent{\bf Theorem~\ref{adjoints}}. {\em Let ${\cal S}$ be a compactly
generated triangulated category, ${\cal T}$ any triangulated category.
Let $F: {\cal S} \rarrow {\cal T}$ be a triangulated functor.  Suppose
$F$ respects coproducts; that is, the natural map
$$
\coprod_{\lambda \in \Lambda} F(s_\lambda) \rarrow
F\left(\coprod_{\lambda \in \Lambda} s_\lambda\right)
$$
is an isomorphism.  Note that although we are not assuming that ${\cal T}$
has coproducts, we are assuming that the object on the right is a coproduct
of the objects on the left.
Then there exists a right adjoint for $F$, namely a functor
$G: {\cal T} \rarrow {\cal S}$, for which there is a natural isomorphism

$$
Hom_{\cal S}\left(x,Gy\right)=Hom_{\cal T}\left(Fx,y\right).
$$ }
\bigskip

\nin
We will see that Grothendieck duality is an immediate consequence. Even
the non--noetherian statements come very cheaply.

Sections~\ref{S1} and \ref{S1.5}  are introductory.  They
give the definitions and elementary
properties of compactly generated categories.  They also discuss why
the categories one naturally wants to consider are examples.  In particular,
the derived category of quasi-coherent sheaves on $X$, and the derived
category of quasi-coherent sheaves of $D$-modules. There is,
however, a technical point. In the context of Brown representability,
it is essential to have triangulated categories with direct sums.
Thus we must work with {\em unbounded} derived categories.
Sections~\ref{S1} and \ref{S1.5} also discuss why
the natural functors one might consider, for instance $R f_*$, respect
coproducts.

It should definitely be noted that nothing in Section~\ref{S1}
is new. Except
for the terminology, the results can certainly all be found in SGA6. There is
also an excellent exposition of them in Thomason's~\cite{TT}. But there are
two reasons for giving a complete and self--contained exposition of these
facts. One is to keep this article self--contained. But more importantly,
both \cite{SGA6} and \cite{TT} are pre--Bousfield. This needs to be made
precise.
As the reader will easily observe by studying the dates, \cite{TT} came
long after Bousfield's \cite{B1}, \cite{B2} and \cite{B3}. But as far as the
author knows, it was not until \cite{BN} that it was realised that Bousfield's
techniques applied to these problems. And \cite{BN} is more recent than
\cite{TT}. Ever since \cite{BN}, the present author has delighted in taking
every  opportunity to point out that Bousfield's techniques make all the old
results much clearer, more general and easier to prove. Section~\ref{S1}
is to be taken in this vein. In fact, perhaps the entire article is little more
than a manifestation of the above.

As for Section~\ref{S1.5}, it redoes the classical theory from the new
perspective offered by Thomason's localisation theorem.
Classically, the proofs that $D(qc/X)$ is compactly generated relied
on finding enough line bundles on $X$. It follows from
Thomason's localisation theorem that one does not need line bundles.
This is explained in Section~\ref{S1.5}. Again, the observation is not
new, it was first made by Thomason in \cite{TT}. Section~\ref{S1.5} is
the only section in this article which assumes familiarity either
with Thomason's article on the subject, or with \cite{STT}. There is,
for the reader's benefit, a summary of the needed results in
Theorem~\ref{Thomason's localisation}. A reader willing to assume that
all his/her schemes have ample line bundles can skip Section~\ref{S1.5},
except for the last page and a half which discuss homotopy colimits.

There is a technical point which perhaps deserves mention. In the literature,
one
frequently considers some other triangulated categories. For instance, the
category $D(qc/X)$ of complexes of quasi--coherent sheaves on $X$ may be
replaced
by $D^{}_{qc}(X)$, the category of complexes of sheaves of ${\cal O}^{}_X$
modules on $X$, with quasi--coherent cohomology. There is a natural functor
$$
F:D(qc/X)\longrightarrow D^{}_{qc}(X).
$$
It turns out that this functor is an equivalence
of categories if $X$ is
quasi--compact and separated; see \cite{BN}, Corollary 5.5. A similar statement
holds
for the derived categories of $D$--modules; I leave it to the reader to state
and
prove the analogue. We will make no explicit use of this fact. It is only
mentioned to
comfort the experts who might be concerned about such things.

A similar statement also is valid for maps. Given a morphism $f:X \la Y$, there
are induced maps $Rf_*:D(qc/X)\la D(qc/Y)$ and
$Rf_*:D(X)\la D(Y)$. It is comforting to know that they agree;
there is a commutative diagram
$$
\begin{array}{ccc}
D(qc/X)  &  \la   &  D(qc/Y) \\
\downarrow&       &\downarrow\\
D(X)  &  \la   &  D(Y).
\end{array}
$$
It is possible to give a proof based on homotopy theoretic ideas
as in \cite{BN}, but for a written proof we refer the reader to
\cite{L1}, Proposition~3.9.2.

Sections \ref{S2} and \ref{S3}
contain the proof of the Brown Representability Theorem,
Theorem~\ref{adjoints}.
The proof itself
is very short and simple, and although not very different from Brown's
proof in \cite{B2}, we include it for the reader's convenience.

The last two sections, Sections \ref{S4} and \ref{S5},
treat the same problems that are addressed
by Verdier in \cite{V}. What is different here is (1) that we work with
unbounded derived
categories, and (2) that the argument is entirely based on the behaviour of
coproducts. In Section~\ref{S4} we ask when does the right adjoint $f^!$ of
$R f_*: D(qc/X) \rarrow D(qc/Y)$ respect coproducts?  It turns out
that one can give a simple, satisfactory criterion, and what makes the
question ``natural'' is that $f^!$ respects coproducts precisely when
$$
f^!(y) \simeq Lf^*(y) \otimes f^! {\cal O}_Y.
$$
In other words, $f^!$ respects coproducts when there is a dualizing complex
$f^!{\cal O}_Y$, and $f^!$ is given as tensor product with $f^!{\cal O}_Y$.
Note that we will always write $\otimes$ for the left derived
functor $^L\otimes$.

The final question we address is when the functor $f^!$ localizes well;
that is, when it gives an isomorphism in $D(qc/Y)$, as in Grothendieck's
original theorem.  Although our treatment here is not complete, we give a
useful sufficient criterion. The reader should note that
Theorem 2 in Verdier~\cite{V},
page 394 seems better than the result we get here; but this is partly a
delusion.
Verdier's theorem is about $D^+_{qc}(X)$ whereas Proposition~\ref{the
vanishing}
deals with $D^{}_{qc}(X)$. The difference between bounded and unbounded derived
categories is crucial here. Without boundedness hypotheses
Theorem 2 in Verdier~\cite{V} fails. Another way to say this is that certain
functors that come up in the proof commute with coproducts in $D^+_{qc}(X)$ but
not with coproducts in $D^{}_{qc}(X)$. See the proof of
Proposition~\ref{the vanishing, bounded} for a proof of
Verdier's result based on coproducts
(in the non--noetherian case too).
See Example~\ref{counterexample} for a counterexample showing
that without some conditions, one cannot expect Verdier's theorem
to hold in the unbounded derived category.

It should emphasised that here also, the topological techniques
work without any hypothesis that the schemes be noetherian. We
recover, therefore, Lipman's results.

The point of this article is that the Grothendieck duality is very
easy to prove by homotopy theoretic techniques. It is an immediate
consequence of Brown representability. But the reader is not
assumed familiar with homotopy theory; hence the first few sections,
in which we attempt to give something like a self--contained treatment.
The reader is not assumed familiar with much of the literature
on the subject, especially not with previous articles by the present
author.  What is assumed is some familiarity with derived categories.
Since we will be working almost entirely with unbounded derived categories,
the reader should be familiar with the work of Spaltenstein \cite{S} on
extending $Rf_*$ and tensor products to unbounded complexes.
There is a brief account of Spaltenstein's results
in Sections 1 and 2 of \cite{BN};
this account is recommended because it uses the notations and terminology
of the topologists, which we will also follow here.

Finally, I wish to thank Morava for directing me to the problem.  I read
{\em Residues and Duality} many years ago, before I learned any topology.
But Morava kept suggesting that there is more there than meets the eye.
There is a topological analogue one would like to understand.  On rereading
{\em Residues and Duality}, it became clear to me that even the algebraic
geometry could be understood better.
I also want to thank Alonso, Jeremias, Kuhn
and Lipman for helpful conversations and useful comments. Lipman
was especially helpful, reading earlier versions of this manuscript,
pointing out gaps and making many suggestions. The referee also pointed out
helpful simplifications and corrections.

The $D$-module problem of directly constructing the right
adjoint of $R f_+$ was posed by Nick Katz 9 years ago.  By a
happy accident, I am now able to solve it.

\section{Preliminaries, approached classically}
\label{S1}

Let ${\cal T}$ be a triangulated category.  There are several hypotheses
one likes to place on ${\cal T}$ in order to work with it.

\dfn{sums exist}  The category ${\cal T}$ is said to
{\em contain small coproducts} if, for any small set $\Lambda$ and any
collection $\{t_\lambda, \: \lambda \in \Lambda\}$ of objects $t_\lambda
\in \mbox{Ob}({\cal T})$ indexed by $\Lambda$, the categorical coproduct
$$
\displaystyle\coprod_{\lambda \in \Lambda} t_\lambda
$$
exists in ${\cal T}$.\edfn

\rmk{R1.comm} It turns out that one can prove two things:
\begin{description}
\sthm{R1.1.1} The suspension functor commutes with coproducts; that is
the natural map
$$
\coprod_{\lambda \in \Lambda} \Sigma t_\lambda \la
\Sigma\left\{\coprod_{\lambda \in \Lambda} t_\lambda\right\}
$$
is an isomorphism.\esthm
\sthm{R1.1.2} The coproduct of any set of triangles is a triangle.
\esthm
\end{description}
The proof of this may be found in the yet--to--be--completed
\cite{NV}. If the reader is unhappy with this, just
modify Definition~\ref{sums exist}
so that a category $\ct$ containing
small coproducts is assumed to
satisfy \ref{R1.1.1} and \ref{R1.1.2}.
\ermk

\exm{sums in $D(qc/X)$} Let $X$ be a scheme.  Let
$D(qc/X)$ be the derived category of chain complexes of quasi-coherent
sheaves over $X$.

Suppose $\{ x_\lambda |\: \lambda \in \Lambda \}$ is a set of objects
of $D(qc/X)$.  That is, each $x_\lambda$ is a chain complex
$$
\rarrow x_\lambda^{n-1} \stackrel{\partial}{\rarrow} x_\lambda^n
\stackrel{\partial}{\rarrow} x_\lambda^{n+1} \rarrow .
$$
Then it is nearly trivial that the chain complex
$$
\longrightarrow \displaystyle\coprod_{\lambda \in \Lambda}
x_\lambda^{n-1} \stackrel{\coprod\partial}{\longrightarrow}
\displaystyle\coprod_{\lambda \in \Lambda} x_\lambda^n
\stackrel{\coprod\partial}{\longrightarrow}
\displaystyle\coprod_{\lambda \in \Lambda}
x_\lambda^{n+1} \longrightarrow
$$
is the coproduct of the $x_\lambda$ in the category $D(qc/X)$.  It is also
nearly trivial that coproducts of triangles are triangles.  For a proof, see
\cite{BN}, Corollary 1.7.\eexm

\lem{$Rf_*$ respects sums}  Let $X$ be a quasi-compact,
separated scheme and $Y$ be a scheme.  Let $f: \: X \rarrow Y$ be a
separated morphism.  Let $R f_*: \:
D(qc/X) \rarrow D(qc/Y)$ be the direct image functor.  Then the natural
map
$$
\displaystyle\coprod_{\lambda \in \Lambda} R f_* x_\lambda
\rarrow R f_*\left\{\displaystyle\coprod_{\lambda \in \Lambda} \:
x_\lambda\right\}
$$
is an isomorphism; that is, $Rf_*$ respects coproducts.\elem

\noindent{\bf Proof.}  The question being local in $Y$, we may assume $Y$
affine.  Since $X$ is quasi-compact, it may be covered by finitely many
open affines: $X = \displaystyle\bigcup^n_{i=1} U_i$, with
$U_i$ affine.  We will prove the lemma by induction on the number $n$ of
open affines.

If $n = 1$, then $X = U_1$ is affine.  Thus the map
$$
\mbox{Spec}(S) = X \rarrow Y = \mbox{Spec}(R)
$$
corresponds to a ring homomorphism $R \rarrow S$.  The category $D(qc/X)$
is just $D(S)$, the derived category of chain complexes of $S$--modules.  The
functor $Rf_*: D(S) \rarrow D(R)$ just takes a chain complex of $S$--modules
and views it as a chain complex of $R$--modules.  This clearly preserves
coproducts.

If $n > 1,$ let $U = U_1, \ V = \displaystyle\bigcup^{n}_{i=2} U_i$.
Then $U \cap V = \displaystyle\bigcup^{n}_{i=2} (U_1 \cap U_i)$, and
$U_1 \cap U_i$ is affine because $X$ is separated.  Thus both $V$ and
$U \cap V$ are unions of $n-1$ affines.  By induction, the theorem holds
for the maps $f_U: U \rarrow Y, \ f_V:V \rarrow Y$ and
$f_{U \cap V}: U \cap V \rarrow Y$.

Let $i_U: U \hookrightarrow X, \ i_V: V \hookrightarrow X$, and
$i_{U \cap V}: U \cap V \hookrightarrow X$ be the open immersions.  Then any
object $x$ of $D(qc/X)$ admits a triangle

\begin{picture}(400,100)
\put(150,10) {$R(i^{}_{U \cap V})_* i_{U \cap V}^* (x)$}
\put(180,70) {$\longrightarrow$}
\put(120,70) {$x$}
\put(220,70) {$R(i^{}_U)_* i^*_U x \oplus R(i^{}_V)_* i^*_V(x)$}
\put(125,34) {$\scriptstyle{(1)}$}
\put(135,40) {$\nwarrow$}
\put(230,40) {$\swarrow$}
\end{picture}

\nin
Thus we deduce a triangle

\begin{picture}(400,100)
\put(140,10) {$Rf_* R(i^{}_{U \cap V})_* i_{U \cap V}^*(x)$}
\put(180,70) {$\longrightarrow$}
\put(110,70) {$Rf_*(x)$}
\put(230,70) {$Rf_*R(i^{}_U)_*i^*_U(x) \oplus Rf_* R(i^{}_V)_*i^*_V(x)$}
\put(125,34) {$\scriptstyle{(1)}$}
\put(135,40) {$\nwarrow$}
\put(230,40) {$\swarrow$}
\end{picture}

\noindent
But now
\begin{eqnarray*}
Rf_* R(i^{}_U)_* & = & R(f^{}_U)_* \\
Rf_* R(i^{}_V)_* & = & R(f^{}_V)_* \\
\mbox{and} \rule{1cm}{0cm} Rf_* R(i^{}_{U \cap V}) & = & R(f^{}_{U \cap V})_*
\rule{1cm}{0cm}
\end{eqnarray*}
all commute with coproducts by the induction hypothesis.  The functors
$i^*_U, \ i^*_V$, and $i^*_{U \cap V}$ commute with coproducts because they
have right adjoints.  Therefore, in the morphism of triangles
$$
\begin{array}{ccccc}
\displaystyle\coprod_{\lambda \in \Lambda} Rf_*(x_\lambda)
& \rightarrow & \displaystyle\coprod_{\lambda \in \Lambda}
[R(f^{}_U)_* i^*_U(x_\lambda) \oplus R(f^{}_V)_* i^*_V(x_\lambda)]
& \rightarrow & \displaystyle\coprod_{\lambda \in \Lambda}
R(f^{}_{U \cap V})_* i^*_{U \cap V} (x_\lambda) \\
\alpha \downarrow && \beta \downarrow && \gamma \downarrow \\
Rf_*\left(\displaystyle\coprod_{\lambda \in \Lambda} x_\lambda\right)
& \rightarrow & R(f^{}_U)_* i^*_U \left(
\displaystyle\coprod_{\lambda \in \Lambda} x_\lambda\right) \oplus R(f^{}_V)_*
i^*_V \left(\displaystyle\coprod_{\lambda \in \Lambda} x_\lambda\right)
& \rightarrow & R(f^{}_{U \cap V})_* i^*_{U \cap V}
\left(\displaystyle\coprod_{\lambda \in \Lambda} x_\lambda\right)
\end{array}
$$
the maps $\beta$ and $\gamma$ are isomorphisms; hence so is $\alpha$.
\hfill $\Box$
\bigskip

\cor{sums} If in Lemma~\ref{$Rf_*$ respects sums}
we take $Y$ to be the scheme Spec$(\boldZ)$, then
$Rf_*$ is just the derived functor of the global section
functor.  We deduce that if $X$ is a quasi-compact,
separated scheme then the functors $H^i$
$$
H^i(X,-): D(qc/X) \rarrow \{\mbox{abelian groups}\}
$$
respect coproducts.\ecor

\dfn{compact objects} An object $c$ of ${\cal T}$
is called {\em compact} if, for any coproduct of objects of ${\cal T}$
$$
\mbox{Hom}_{\cal T} \left(c, \displaystyle\coprod_{\lambda \in \Lambda}
t_\lambda \right) = \displaystyle\coprod_{\lambda \in \Lambda}
\mbox{Hom}_{\cal T} (c, t_\lambda).
$$
\edfn
\medskip

\noindent{\bf Observation.}  The suspension of a compact object is compact.

\dfn{compactly generated}  The triangulated category
${\cal T}$ is called {\em compactly generated} if ${\cal T}$ contains small
coproducts,
and there exists a small set $T$ of compact objects of ${\cal T}$, such that
$$
\mbox{Hom} (T,x) = 0 \Rightarrow x = 0.
$$
In other words, if $x$ is an object of ${\cal T}$, and for every $t \in T$,
Hom$(t,x) = 0$, then $x$ must be the zero object.\edfn

\dfn{generating set}  If ${\cal T}$ is a
compactly generated triangulated category, then a set $T$ of compact objects
of ${\cal T}$ is called a {\em generating set} if
\begin{itemize}
\item[(1)] Hom$(T,x) = 0 \Rightarrow x = 0;$
\item[(2)] $T$ is closed under suspension; $T = \Sigma T$.
\end{itemize}\edfn

\rmk{2 is easy} Let $T$ be any set of objects in
${\cal T}$ as in Definition~\ref{compactly generated}. Then
$$
\displaystyle\bigcup_{i \in \boldZ} \Sigma^i T
$$
is a generating set as in Definition~\ref{generating set}.  (1) holds because
it holds for $T$; and clearly the set
$$
\displaystyle\bigcup_{i \in \boldZ} \Sigma^iT
$$
is stable under suspension.  It needs only be remarked that any suspension of
a compact object is compact, and hence the set consists only of compact
objects of ${\cal T}$.\ermk

\exm{quasi-projective variety} Let $X$ be a quasi-compact,
separated scheme.  Let ${\cal T} = D(qc/X)$ be the category of chain
complexes of quasi-coherent sheaves on $X$.  By
Example~\ref{sums in $D(qc/X)$} we know that
${\cal T}$ contains small coproducts.  Now let
${\cal L}$ be any line bundle on $X$.  View ${\cal L}$ as an object of
$D(qc/X)$; it is the complex of sheaves which is just ${\cal L}$ in degree
0. Then
$$
\mbox{Hom}\left({\cal L}, \displaystyle\coprod_{\lambda \in \Lambda}
t_\lambda\right) = H^0
\left({\cal L}^{-1} \otimes \displaystyle\coprod_{\lambda \in \Lambda}
t_\lambda\right).
$$
Now tensor product respect coproducts, since it has a right adjoint.
The functor $H^0$ respects coproducts by Corollary~\ref{sums}. It follows
that ${\cal L}$ is a compact object of ${\cal T} = D(qc/X)$.

Now suppose ${\cal L}$ is an ample line bundle.  For any $m \in \boldZ, \
{\cal L}^m$ is compact.  For any $n \in \boldZ, \ \Sigma^n {\cal L}^m$ is
also compact. Let
$$
T = \{\Sigma^n {\cal L}^m| \: m,n \in \boldZ\}.
$$
I assert that $T$ is a generating set for ${\cal T}$. Suppose $x \neq 0$
is an object of ${\cal T}$. Then it has some non-trivial sheaf cohomology;
for some $n, \ {\cal H}^{-n}(x) \neq 0$.  But ${\cal H}^{-n}(x)$ is
a quasi-coherent sheaf on $X$, and because ${\cal L}$ is ample,
${\cal L}^t \otimes {\cal H}^{-n}(x)$ has non-zero global sections
for some $t >> 0$.  If $x$ is the complex
$$
\rarrow x^{-n-1} \stackrel{\partial}{\rarrow} x^{-n} \stackrel{\partial}
{\rarrow} x^{-n+1} \rarrow,
$$
then there is a surjective map of quasi-coherent sheaves on $X$
$$
\mbox{ker}(x^{-n} \rarrow x^{-n+1}) \rarrow {\cal H}^{-n}(x),
$$
and it follows that, choosing $t$ large enough, we can find a class
$$
s \in H^0({\cal L}^t \otimes x^{-n})
$$
which maps to 0 in $H^0({\cal L}^t \otimes x^{-n+1})$, and whose image in
$H^0({\cal H}^{-n}({\cal L}^t \otimes x))$ is non-zero.  This
immediately gives us a non-zero map
$$
\Sigma^n {\cal L}^{-t} \rarrow x.
$$
If all such maps vanish, so must $x$.\eexm

\exm{families of line bundles} In Example~\ref{quasi-projective variety}
we can replace a single ample line bundle by
a family; if $X$ admits a family of line bundles
$\{{\cal L}_\alpha | \: \alpha \in A\}$ which is jointly ample, then
$$
T = \{\Sigma^n {\cal L}^m_\alpha | \: m,n \in \boldZ, \: \alpha \in A\}
$$
will do for a generating set of compact objects. See
\cite{SGA6}, page 171.\eexm

\exm{$D$-modules}  Let $X$ be a smooth quasi-compact
variety of finite type over a field $k$ of characteristic zero.
Let ${\cal T} = D\left(\frac{qc\: D\mbox{-modules}}{X}\right)$ be
the derived category of chain complexes of quasi-coherent $D$-modules over
$X$.  Once again, it is trivial that ${\cal T}$ contains small coproducts.

Because $X$ is smooth, one knows by a trick of Kleiman that there is an
ample family of line bundles; cover $X$ by open affines
$X = \displaystyle \bigcup_{i\in I} U_i$
and let ${\cal L}_i$ be the line bundle ${\cal O}(D_i)$, where $D_i$ is
the divisor $X-U_i$.  Then the ${\cal L}_i$'s form an ample family.

Now observe
$$
\mbox{Hom}_{D\left(\frac{qc\: D\mbox{-modules}}{X}\right)}
\left({\cal D}^{}_X \otimes_{{\cal O}^{}_X} {\cal L}_i^m, x\right)
= \mbox{Hom}_{D\left(qc/X\right)} \left({\cal L}_i^m,x\right).
$$

\nin
 From this it follows that
$\Sigma^n\left({\cal D}^{}_X \otimes_{{\cal O}^{}_X} {\cal L}_i^m\right)$
are compact for all $i, m$ and $n$, and the set
$$
T = \left\{\Sigma^n
\left({\cal D}^{}_X \otimes_{{\cal O}^{}_X} {\cal L}_i^m\right)\left| \:
i\in I,\quad m,n \in \boldZ\right.\right\}
$$
is a generating set.\eexm

\exm{other compacts in $D(qc/X)$}
Let $X$ be a quasi-compact, separated scheme.  In
Example~\ref{quasi-projective variety}
we proved that any line bundle ${\cal L}$ on $X$,
viewed as an object of $D\left(qc/X\right)$, is compact.
Given ample families of line bundles on $X$, we used this to construct a
generating set.

Let us now observe that if $c \in D\left(qc/X\right)$ is any perfect
complex on $X$, it is compact.  Recall that a complex $c$ is perfect if,
locally
on $X$, it is isomorphic to a bounded complex of finitely generated,
projective ${\cal O}^{}_X$-modules.\eexm

\noindent{\bf Proof.}  Let $\displaystyle\coprod_{\lambda \in \Lambda} \:
x_\lambda$ be a coproduct in $D\left(qc/X\right)$.
Let
$$
{\cal RH}om\left(c, \displaystyle\coprod_{\lambda \in \Lambda}
x_\lambda\right)
$$
be the sheaf ${\cal RH}om$; it is an object of $D\left(qc/X\right)$.
I assert that the natural map in $D\left(qc/X\right)$
$$
\phi_c: \: \coprod_{\lambda \in \Lambda} {\cal RH}om (c, x_\lambda)
 \rightarrow {\cal RH}om\left(c, \coprod_{\lambda \in \Lambda} x_\lambda\right)
$$
is an isomorphism whenever $c$ is a perfect complex.
The problem is local in $X$; we may therefore assume
that $X$ is affine, and $c$ is a bounded complex of finitely generated
projective ${\cal O}_X$-modules.

If $c$ is $\Sigma^m{\cal O}_X$, then $\phi_c$ is clearly an isomorphism.
If $c$ is $\Sigma^m{\cal O}^n_X$, a finite direct sum of $\Sigma^m
{\cal O}_X$'s, then $\phi_c$ is also an isomorphism.  If $c = c' \oplus
c''$ and $\phi_c$ is an isomorphism, so are $\phi_{c'}$ and $\phi_{c''}$.
Hence $\phi_c$ is an isomorphism whenever $c$ is a suspension of a
finitely generated projective module.

Now if we have a triangle

\begin{picture}(400,100)
\put(182,10) {$c''$}
\put(180,70) {$\longrightarrow$}
\put(135,70) {$c$}
\put(230,70) {$c'$}
\put(145,34) {$\scriptstyle{(1)}$}
\put(150,40) {$\nwarrow$}
\put(210,40) {$\swarrow$}
\end{picture}

\noindent
and $\phi_{\Sigma^mc'}$ and $\phi_{\Sigma^mc''}$ are isomorphisms for
all $m \in \boldZ$, then it follows from the 5-lemma that $\phi_{\Sigma^mc}$
is an isomorphism for all $m \in \boldZ$. Therefore the full subcategory of
$c$'s such that $\phi_{\Sigma^mc}$ is an isomorphism for all $m
\in \boldZ$ is triangulated and contains the finitely generated projective
${\cal O}_X$-modules.  Hence it contains finite complexes of finitely
generated projectives.

Thus we have proved that for $c$ a perfect complex,
$$
\phi_c : \ \displaystyle\coprod_{\lambda \in \Lambda}
{\cal RH}om (c,x_\lambda) \rightarrow
{\cal RH}om\left(c, \displaystyle\coprod_{\lambda \in \Lambda} x_\lambda\right)
$$
is an isomorphism.  But
\begin{eqnarray*}
\mbox{Hom}\left(c, \displaystyle\coprod_{\lambda \in \Lambda} x_\lambda\right)
& = & H^0\left[{\cal RH}om\left(c, \coprod_{\lambda \in \Lambda} x_\lambda
\right)\right] \\
& = & H^0\left[\coprod_{\lambda \in \Lambda}{\cal RH}om\left(c, x_\lambda
\right)\right] \\
& = & \coprod_{\lambda \in \Lambda}H^0\left[{\cal RH}om\left(c, x_\lambda
\right)\right] \\
& = & \coprod_{\lambda \in \Lambda}\mbox{Hom}\left(c,x_\lambda\right),
\end{eqnarray*}
where the third equality is by Corollary~\ref{sums}, which assures us that
$H^0$
respects coproducts.  \hfill $\Box$
\bigskip

\exm{other compacts in $D(qc D-modules/X)$}
Let $X$ be a quasi-compact, separated smooth scheme of finite type over
a field $k$ of characteristic 0.  Then by Example~\ref{$D$-modules} we know
that objects of the form $\Sigma^n\{{\cal D}_X \otimes_{{\cal O}_{X}}
{\cal L}\}$ are compact in ${\cal T} =
D\left(\frac{qc\: D\mbox{-modules}}{X}\right)$. But
much as in Example~\ref{other compacts in $D(qc/X)$}, it can be shown that
any bounded complex with coherent cohomology is compact.  Note that because
${\cal D}_X$ is locally of finite projective dimension, any coherent sheaf
can locally be replaced by a finite chain complex of finitely generated
projectives.  We leave the details to the reader.\eexm

\section{The approach using Thomason's localisation theorem}
\label{S1.5}

Let $\cal T$ be a triangulated category. In this section, and for the
remainder of the article, we adopt the notation that ${\cal T}^c$
stands for the full subcategory of compact objects in $\cal T$.

Let $X$ be a quasi--compact, separated scheme. In Section~\ref{S1} we
saw how to prove that the category $D(qc/X)$ is compactly generated,
provided $X$ has an ample family of line bundles; see
Example~\ref{families of line bundles}. On the other hand, we
also know that if $X$ is arbitrary (that is, quasi--compact and separated
but not necessarily possessing any line bundles), every perfect
complex on $X$ is compact. In this section we will use Thomason's
localisation theorem to prove that $D(qc/X)$ is compactly generated,
even without line bundles. First, we quote the theorem to which we
will appeal

\thm{Thomason's localisation} Let $\cal S$ be a
compactly generated triangulated category. Let $R$ be a set of
compact objects of $\cal S$ closed under suspension. Let
$\cal R$ be the smallest full subcategory of $\cal S$
containing $R$ and closed
with respect to coproducts and triangles. Let $\cal T$ be the
Verdier quotient ${\cal S}/{\cal R}$. Then we know:

\sthm{Thom1} The category $\cal R$ is compactly generated,
with $R$ as a generating set.\esthm

\sthm{Thom1.5} If $R$ happens to be a generating set for all of
$\cal S$, then ${\cal R}={\cal S}$.\esthm

\sthm{Thom2} If $R\subset\cal R$ is closed under the formation
of triangles and direct summands, then it is all of ${\cal R}^c$.
In any case, ${\cal R}^c={\cal R}\cap{\cal S}^c$.
\esthm

\sthm{Thom3} Suppose $t$ is a compact object of $\cal T$. Then
there is an object $t'\in{\cal T}^c$ and an object $s\in {\cal S}^c$
and an isomorphism in $\cal T$
$$
s\simeq t\oplus t'.
$$
Thus, $t$ might not be the isomorphic in $\cal T$
to a compact object in $\cal S$,
but it is the direct summand of an object isomorphic
in $\cal T$ to a compact object of $\cal S$.
Furthermore, $t'$ may be chosen to be $\Sigma t$, or
any other object whose sum with $t$ is zero in $K_0$.\esthm

\sthm{Thom4} Given an object $s\in {\cal S}^c$, an object
$s'\in\cal S$, and a morphism in $\cal T$ $s\la s'$,
then there is a diagram in $\cal S$
$$
\begin{array}{ccccccccc}
    &     &   \tilde s  &    &    \\
    &\swarrow&          &\searrow& \\
s   &        &          &      & s'
\end{array}
$$
where $\tilde s$ is compact, in the triangle $r\la \tilde s\la s\la \Sigma r$,
the object $r$ is in ${\cal R}^c$, and when we reduce the diagram
to $\cal T$, the composite of the map $\tilde s\la s'$ with
the inverse of  $\tilde s\la s$ is the given map $s\la s'$.\esthm
\ethm

\rmk{The proofs} For triangulated categories $D(qc/X)$ where $X$
is a scheme, the theorem is due to Thomason \cite{TT}. In the
generality in which the theorem is stated, it may be found in \cite{STT}.
Note that in \cite{STT} we assume not only that $\cal S$ is compactly
generated, but that the generating set may be taken to be ${\cal S}^c$;
in particular, we assume ${\cal S}^c$ to be small. The reader will note
that this is inessential to any of the arguments in \cite{STT}. The only
point where the proof uses the hypothesis on the smallness of
${\cal S}^c$ is in showing \ref{Thom3}, and that comes at the very end.
So in any case we know the other properties. And in the proof of \ref{Thom3},
it suffices to know that $\cal S$ is compactly generated, and that by
\ref{Thom1.5} it then follows that if $S$ is a generating set, the category
${\cal R}\subset{\cal S}$ which is the smallest category containing $S$ and
closed with respect to triangles and coproducts is all of $\cal S$.

With this comment, we now give the references for the proofs.
\ref{Thom1.5} really follows from the proof in \cite{STT}. In the
notation there, if $j^*:{\cal S}\la{\cal T}$ is the natural functor
and $j_*:{\cal T}\la{\cal S}$ its right adjoint, then the objects
$j_*j^*x$ are characterised by the fact that $Hom(R,j_*j^*x)=0$.
Since $R$ generates $\cal S$, it follows that for every $x$,
$j_*j^*x=0$. But the identity on $j^*x$ factors as
$$
j^*x\la j^*j_*j^*x\la j^*x
$$
and the middle object is zero. Hence $j^*x=0$. This is true for every $x$,
which means that the category $\cal T$ is zero, and hence ${\cal R}={\cal S}$.

For the remaining statements,
\ref{Thom2} is Lemma~2.2 of \cite{STT}.
\ref{Thom3} goes as follows. The existence of $t'$ is Lemma~2.6 of
\cite{STT}, but as noted above, the proof needs to be modified
slightly to account for the fact that we are only assuming $\cal S$
compactly generated. The fact that $t'$ may be taken to be
$\Sigma t$, or even anything else which is isomorphic to $-t$ in
$K_0$, may be found in the appendix to \cite{STT}.
\ref{Thom4} is Lemma~2.5 of \cite{STT}.
This leaves us with \ref{Thom1}, which is trivial, so
let us give the proof. Suppose $r$ is an object of $\cal R$ such
that $Hom(R,r)=0$. Consider the full subcategory $^\perp r\subset\cal R$,
defined by
$$
^\perp r=\left\{x\in {\cal R}|Hom(\Sigma^n x,r)=0 \hbox{\rm\ for all\ }
 n\in{\Bbb Z}
\right\}.
$$
Clearly, $^\perp r$ is triangulated and closed under coproducts, and
contains $R$. Hence it must be all of $\cal R$, in particular $^\perp r$
contains $r\in\cal R$. Thus $Hom(r,r)=0$, hence $r=0$.
\ermk

\cor{all the compacts are perf} Let $X$ be a quasi--compact, separated
scheme, and suppose we know that $D(qc/X)$ is compactly generated,
and that the generating set consists of some perfect complexes.
Then the category of all perfect complexes on $X$ is nothing other
than  $D(qc/X)^c$.\ecor

\rmk{for now} For now, we only know that $D(qc/X)$ is compactly generated
when there is an ample family of line bundles. So for now
Corollary~\ref{all the compacts are perf} only applies in that
case. But this will change by the end of the section.\ermk

\nin
{\bf Proof of Corollary~\ref{all the compacts are perf}.}\ \
Let ${\cal S}=D(qc/X)$, and let
$R$ be the set of perfect complexes
in Theorem~\ref{Thomason's localisation}.
By hypothesis, $R$ generates $\cal S$, and hence by
\ref{Thom1.5}, ${\cal R}={\cal S}=D(qc/X)$, and $\cal T$ is
trivial.
But $R$ is closed with respect to direct summands
 and triangles. Closure
with respect to triangles is clear. Closure with respect to
direct summands asserts that a direct summand of a perfect
complex is perfect. This is local, so we may assume $X$ affine.
For affine $X$ this is Proposition~3.4 in \cite{BN}.
By \ref{Thom2} we therefore know that $R={\cal R}^c$; every
compact object in $D(qc/X)$ is a perfect complex.\qqed

Now we come to the existence of compacts on a general $X$.

\pro{compacts exist in general} Let $X$ be a quasi--compact,
separated scheme. Then the category $D(qc/X)$ is compactly generated.\epro

\nin
It might be useful to state a lemma that clearly implies
Proposition~\ref{compacts exist in general}, and which we will
prove.

\lem{extending maps, special case} Let $X$ be a quasi--compact,
separated scheme. Let $U\subset X$ be a quasi--compact, open
subscheme. Let $x$ be an arbitrary object of $D(qc/X)$, and
let $u$ be a perfect complex in $D(qc/U)$. Suppose that we are
given a map in $D(qc/U)$ of the form $u\la x$. Then there
exists a perfect complex $u'$ in $D(qc/U)$ so that the map
$$
u\oplus u'\stackrel{\pi_1}\la u\la x
$$
lifts to $D(qc/X)$. There exists a perfect complex $\tilde u$
on $X$, restricting to $u\oplus u'$ on $U$, and a map
$\tilde u\la x$, defined on $X$,
which restricts on $U$ to the given map
$u\oplus u'\stackrel{\pi_1}\la u\la x$.\elem

\nin
{\bf Proof that Lemma~\ref{extending maps, special case}
implies
Proposition~\ref{compacts exist in general}.}\ \
Take $U$ to be an open affine. Then $U$ admits an ample family
of line bundles; the trivial bundle is ample. Thus we already know
that $U$ is compactly generated, and that the compact objects
are the perfect complexes. If $x$ is any object of $D(qc/X)$
and the restriction of $x$ to $U$ is non-zero, there is a
perfect complex $u$ on $U$ and a non--zero map $u\la x$ on $U$.
By Lemma~\ref{extending maps, special case} this map can
be extended to a non--zero map $\tilde u\la x$ on all of $X$,
where $\tilde u$ is perfect. Thus, unless the restriction of $x$
to every open affine $U\subset X$ is zero, there is a non--zero
map from a perfect complex to $x$. But if the restriction of $x$
to every open affine $U\subset X$ vanishes, then $x$ vanishes.\qqed

\nin
{\bf Proof of Lemma~\ref{extending maps, special case}.}\ \
Suppose first that $X$ is affine.
Then $X$ has an ample line bundle; after all, the trivial bundle is
ample. We therefore already know that $D(qc/X)$ is compactly
generated. Furthermore, the compacts are precisely the
perfect complexes. Now let $D_{X-U}(qc/X)\subset D(qc/X)$
be the full subcategory whose objects are complexes
supported on $X-U$. That is,
$$
D_{X-U}(qc/X)=\left\{x\in D(qc/X)|\hbox{the restriction of $x$ to $U$
is acyclic.}\right\}
$$
Lemma 6.1 of \cite{BN}
shows that even $D_{X-U}(qc/X)$ is compactly
generated, in fact, generated by the suspensions of
one compact object in $D(qc/X)$.

In Theorem~\ref{Thomason's localisation}, let $\cal S$ be the category
$D(qc/X)$, and let $R$ be a generating set for
$D_{X-U}(qc/X)$. This makes ${\cal R}=D_{X-U}(qc/X)$.
Now $\cal T$ is easily identified with $D(qc/U)$.
Thomason localisation theorem applies, and we discover
first that by \ref{Thom3} the complex $u\oplus\Sigma u$
may be lifted to a perfect complex $\tilde u$ in
$D(qc/X)$, and then by \ref{Thom4} that the map
$u\oplus \Sigma u\la x$ can be lifted to a map
$\tilde u\la x$ on all of $X$, possibly after changing
the choice of $\tilde u$ lifting $u\oplus \Sigma u$.

Suppose next that $X=U\cup W$, where $W$ is affine.
We know by the above that the restriction of the
map $u\la x$ can be extended from $U\cap W$ to all of $W$.
Precisely, there is a perfect complex $\tilde u$ on
$W$ and a map $\tilde u\la x$ of complexes on $W$,
so that the restriction to $U\cap W$ is isomorphic
to the map $u\oplus\Sigma u\la x$. Thus, if $j^{}_W:W\la X$,
$j^{}_U:U\la X$ and $j^{}_{U\cap W}:U\cap W\la X$
are the open immersions, we have an isomorphism
on $U\cap W$ of $\{j^{}_{U\cap W}\}^*\{u\oplus\Sigma u\}$ with
$\{j^{}_{U\cap W}\}^*\tilde u$. We are given the maps $\beta$
and $\gamma$, which we can complete to a morphism
of triangles of complexes on $U\cup W=X$
$$
\begin{array}{ccccc}
\hat u
& \rightarrow & \{j^{}_U\}_*\{u\oplus\Sigma u\}\oplus\{j^{}_W\}_*\tilde u
& \rightarrow &  \{j^{}_{U\cap W}\}_*\{j^{}_{U\cap W}\}^*\{u\oplus\Sigma u\}\\
\alpha \downarrow && \beta \downarrow && \gamma \downarrow \\
x
& \rightarrow & \{j^{}_U\}_*\{j^{}_{U}\}^*x\oplus\{j^{}_W\}_*\{j^{}_{W}\}^*x
& \rightarrow & \{j^{}_{U\cap W}\}_*\{j^{}_{U\cap W}\}^*x
\end{array}
$$
and it is easy to check that $\hat u$ is perfect and the map $\hat u\la x$,
defined on all of $X$, is just a lifting of $u\oplus\Sigma u\la x$
from $U$ to all of $X$.

Since $X$ is quasi--compact, it can be covered by finitely many
open affines, and in finitely many steps of extending from
$U$ to $U\cup W$, with $W$ affine, we extend to all of $X$.\qqed

Theorem~\ref{Thomason's localisation}
can also be used to construct
compactly generated categories. The point
being that given a compactly generated category $\cal S$ and
a set of compact objects $R$ in it, the categories $\cal R$
and $\cal T$ are compactly generated.

\exm{holo} Let $X$ be a smooth, quasi-compact, separated
scheme of finite type over a field $k$ of characteristic 0.  Let ${\cal S}
=D\left(\frac{qc\: D\mbox{-modules}}{X}\right) $ be the derived
category of chain complexes
of quasi-coherent $D$-modules on $X$.

Let $R \subset {\cal S}$ be the set
$$
R = \left\{x \in \mbox{Ob}({\cal S})\left| \begin{array}{l}
{\cal H}^i(x) = 0 \mbox{ for all but finitely many } i, \mbox{ and} \\
{\cal H}^i(x) \mbox{ is holonomic for all } i. \end{array} \right.\right\}
$$
Because holonomic modules are coherent, it follows from
Example~\ref{other compacts in $D(qc D-modules/X)$}
that $R$ consists of compact
objects of ${\cal T}$.  Let ${\cal R}$ be as
above.  Then ${\cal R}$
is a compactly generated triangulated
category, with $R$ for a generating set.  We will call this ${\cal R}$ by the
name
$D(\mbox{holo}/X)$.\eexm

The key tool one uses is the homotopy colimit.  Let ${\cal T}$ be a
triangulated category containing small coproducts.  Let
$$
X_0 \rarrow X_1 \rarrow X_2 \rarrow \cdots
$$
be a sequence of objects and morphisms in ${\cal T}$.  Then
$\begin{array}[t]{c} \mbox{hocolim} \\ \vspace{-.75cm} \\
\longrightarrow \end{array} X_i$
is by definition the third edge of the triangle
$$
\begin{array}{rcccc}
\displaystyle\coprod_i X_i & &
\begin{array}{c} 1{\rm -shift} \\ \vspace{-.75cm} \\ \longrightarrow
\end{array}
&& \displaystyle\coprod_i X_i \\
& {\scriptstyle{(1)}} \nwarrow && \swarrow & \\
&& \begin{array}{c} {\rm hocolim} \\ \vspace{-.75cm} \\ \longrightarrow
\end{array} X_i &&
\end{array}
$$
\bigskip

\lem{maps to colimits}  Suppose $c$ is a compact
object of ${\cal T}$, and suppose
$$
X_0 \rarrow X_1 \rarrow X_2 \rarrow \cdots
$$
is a sequence of objects and morphisms in ${\cal T}$.  Suppose ${\cal T}$
admits small coproducts.  Then
$$
\mbox{Hom}\left(c, \begin{array}[t]{c} {\rm hocolim} \\ \vspace{-.75cm} \\
\longrightarrow \end{array} X_i \right) = \displaystyle\lim_{\rarrow}
\mbox{ Hom}(c,X_i).
$$\elem

\noindent{\bf Proof.}  Consider the triangle
$$
\begin{array}{rcccc}
\displaystyle\coprod_i X_i & &
\begin{array}{c} 1{\rm -shift} \\ \vspace{-.75cm} \\ \longrightarrow
\end{array}
&& \displaystyle\coprod_i X_i \\
& {\scriptstyle{(1)}} \nwarrow && \swarrow & \\
&& \begin{array}{c} {\rm hocolim} \\ \vspace{-.75cm} \\ \longrightarrow
\end{array} X_i &&
\end{array}
$$
Applying the homological functor Hom$(c, -)$ we get a long exact
sequence.  In particular
$$
\begin{array}{ccccc}
\mbox{Hom}\left(c, \displaystyle\coprod_i X_i\right)
& \stackrel{\gamma}{\longrightarrow}
& \mbox{Hom}\left(c, \mbox{hocolim} X_i\right)
& \longrightarrow
& \mbox{Hom}\left(c, \displaystyle\coprod_i \Sigma X_i\right) \\
&&&& \rule{0cm}{.01cm} \\
&&&& \downarrow 1\mbox{-shift} \\
&&&& \rule{0cm}{.01cm} \\
&&&& \mbox{Hom}\left(c, \displaystyle\coprod_i \Sigma X_i\right)
\end{array}
$$
is exact.  But $c$ is compact, and hence in the following commutative square
the columns are isomorphisms
$$
\begin{array}{ccc}
\displaystyle\coprod_i \mbox{Hom}\left(c, \Sigma X_i\right)
& \stackrel{1\mbox{-shift}}{\longrightarrow}
& \displaystyle\coprod_i \mbox{Hom}\left(c, \Sigma X_i\right) \\
&& \rule{0cm}{.01cm} \\
|\downarrow \wr && |\downarrow \wr \\
&& \rule{0cm}{.01cm} \\
\mbox{Hom}\left(c, \displaystyle\coprod_i\Sigma X_i\right)
& \stackrel{1\mbox{-shift}}{\longrightarrow}
& \mbox{Hom}\left(c, \displaystyle\coprod_i\Sigma X_i\right).
\end{array}
$$
But the top row is clearly injective. Hence the bottom row is injective,
and we deduce that $\gamma$ is surjective.

We now have a commutative diagram
$$
\begin{array}{ccccccc}
\displaystyle\coprod_i\mbox{Hom}\left(c, X_i\right)
& \stackrel{1\mbox{-shift}}{\longrightarrow}
& \displaystyle\coprod_i \mbox{Hom}\left(c, X_i\right) &&&& \\
&&&& \rule{0cm}{.01cm} &&\\
|\downarrow \wr && |\downarrow \wr && &&\\
&&&& \rule{0cm}{.01cm} &&\\
\mbox{Hom}\left(c, \displaystyle\coprod_iX_i\right)
& \stackrel{1\mbox{-shift}}{\longrightarrow}
& \mbox{Hom}\left(c, \displaystyle\coprod_i X_i\right)
& \stackrel{\gamma}{\longrightarrow} &
\mbox{Hom}\left(c, \begin{array}{c}
\mbox{hocolim} \\ \vspace{-.75cm} \\ \longrightarrow
\end{array} X_i \right) &\longrightarrow&0\\
\end{array}
$$
where the bottom row is exact.  The top row identifies
Hom$\left(c, \begin{array}{c}
\mbox{hocolim} \\ \vspace{-.75cm} \\ \longrightarrow
\end{array} X_i \right)$ as
$\displaystyle\lim_\rightarrow $ Hom$(c,X_i)$. \hfill  $\Box$

\section{Brown Representability}
\label{S2}

In this section, we will prove:
\bigskip

\thm{Brown representability}  Let ${\cal T}$ be a
compactly generated triangulated category.  Let $H: {\cal T}^{op} \rarrow
Ab$ be a homological functor.  That is, $H$ is contravariant and takes
triangles to long exact sequences.  Suppose the natural map
$$
H\left(\displaystyle\coprod_{\lambda \in \Lambda} t_\lambda\right)
\rarrow \displaystyle\prod_{\lambda \in \Lambda} H(t_\lambda)
$$
is an isomorphism for all small coproducts in ${\cal T}$.  Then $H$ is
representable.\ethm

\noindent{\bf Proof.}  Let $T$ be a generating set for ${\cal T}$.  Let
$U_0$ be defined as
$$
U_0 = \displaystyle\bigcup_{t \in T} H(t).
$$
Thus elements of $U_0$ can be thought of as pairs $(\alpha,t)$ with $\alpha
\in H(t)$.  Put
$$
X_0 = \displaystyle\coprod_{(\alpha,t) \in U_{0}} t.
$$
Then
$$
H(X_0) = \displaystyle\prod_{(\alpha,t) \in U_{0}} H(t),
$$
and there is an obvious element in $H(X_0)$, namely the element which is
$\alpha \in H(t)$ for $(\alpha,t) \in U_0$. Call this element $\alpha_0
\in H(X_0)$.  The construction is such that if $t \rarrow X_0$ is the
inclusion of $t$ into $X_0 = \displaystyle\coprod_{(\alpha,t) \in U_{0}} t$
corresponding to $(\alpha, t) \in U_0$, then the induced map $H(X_0) \rarrow
H(t)$
takes $\alpha_0 \in H(X_0)$ to $\alpha \in H(t)$.

To give an object $X_0$ and an element $\alpha_0 \in H(X_0)$ is by Yoneda's
lemma the same as giving a natural transformation
$$
\phi_0: \mbox{ Hom}(-, X_0) \rarrow H,
$$
and what we have seen is precisely that
$$
\phi_0(t): \mbox{ Hom}(t,X_0) \rarrow H(t)
$$
is surjective, for all $t \in T$.

Suppose that for some $i \geq 0$ we have defined an object $X_i$ of ${\cal T}$,
and a natural transformation
$$
\mbox{Hom}(-,X_i) \rarrow H.
$$
Define $U_{i+1}$ by
$$
U_{i+1} = \displaystyle\bigcup_{t \in T} \ker\left\{
\mbox{Hom}(t,X_i) \rarrow H(t)\right\}.
$$
An element of $U_{i+1}$ can be thought of as a pair $(f,t)$ where
$t \in T$ and $f: t \rarrow X_i$ is a morphism. Put
$$
K_{i+1} = \displaystyle\coprod_{(f,t) \in U_{i+1}} t,
$$
and let $K_{i+1} \rarrow X_i$ be the map which is $f$ on the factor $t$
corresponding to the pair $(f,t)$.  Let $X_{i+1}$ be given by the
triangle
$$
\begin{array}{rcccc}
K_{i+1} && \longrightarrow && X_i \\
& {\scriptstyle{(1)}} \nwarrow && \swarrow & \\
&& X_{i+1} &&
\end{array}
$$
We have a map Hom$(-,X_i) \rarrow H$, which by Yoneda's lemma
corresponds to an element $\alpha_i \in H(X_i)$.  Under the map
\begin{eqnarray*}
H(X_i) \rarrow H(K_{i+1}) & = & H\left(\displaystyle\coprod_{(f,t) \in
U_{i+1}} t\right) \\
& = & \displaystyle\prod_{(f,t) \in U_{i+1}} H(t)
\end{eqnarray*}
the element $\alpha_i \in H(X_i)$ maps to zero; the $f: t \rarrow X_i$ were
chosen
so that the induced map Hom$(t,X_i) \rarrow H(t)$ vanishes.  But $H$ is a
homological functor; the exact sequence
$$
H(X_{i+1}) \stackrel{k}{\rarrow} H(X_i) \stackrel{j}{\rarrow}
H(K_{i+1})
$$
coupled with the fact that $j(\alpha_i) =0$, guarantees that there exists
$\alpha_{i+1} \in H(X_{i+1})$ with $k(\alpha_{i+1}) = \alpha_i$.  Choose such
an
$\alpha_{i+1}$.  There is a corresponding natural transformation
$$
\mbox{Hom}(-,X_{i+1}) \rarrow H
$$
rendering commutative the triangle

\begin{picture}(400,100)
\put(110,10) {$\mbox{Hom}(-,X_{i+1})$}
\put(200,10) {$\longrightarrow$}
\put(260,10) {$H$}
\put(180,70) {$\mbox{Hom}(-,X_i)$}
\put(170,40) {$\swarrow$}
\put(240,40) {$\searrow$}
\end{picture}

\nin
Let $X = \begin{array}{c} \mbox{hocolim} \\ \vspace{-.75cm} \\
\longrightarrow \end{array} X_i$.  I assert:

\begin{itemize}
\item[(1)] There is a natural transformation
Hom$(-,X) \rarrow H$ rendering commutative the triangles

\begin{picture}(400,100)
\put(120,10) {$\mbox{Hom}(-,X)$}
\put(200,10) {$\longrightarrow$}
\put(260,10) {$H$}
\put(180,70) {$\mbox{Hom}(-,X_i)$}
\put(170,40) {$\swarrow$}
\put(240,40) {$\searrow$}
\end{picture}

\noindent for every $i$.

\item[(2)] The natural transformation Hom$(-,X) \rarrow H$ is an isomorphism.
\end{itemize}

\noindent{\bf Proof of (1).}  consider the triangle
$$
\begin{array}{rcccc}
\displaystyle\coprod_i X_i && \stackrel{1\mbox{-shift}}{\longrightarrow}
&& \displaystyle\coprod_i X_i  \\
& {\scriptstyle{(1)}} \nwarrow && \swarrow & \\
&& \stackrel{\mbox{hocolim}}{\longrightarrow} X_i = X
&& \end{array}
$$
Applying the cohomological functor $H$, we get an exact sequence
$$
\begin{array}{ccccc}
H(X) & \longrightarrow & H\left(\displaystyle\coprod_i X_i \right) &
\stackrel{1\mbox{-shift}}{\longrightarrow} &
H\left(\displaystyle\coprod_i X_i \right) \\
&&&& \rule{0cm}{.01cm} \\
&& \| && \| \\
&&&& \rule{0cm}{.01cm} \\
&& \displaystyle\prod_i H(X_i)
& \stackrel{1\mbox{-shift}}{\longrightarrow} &
\displaystyle\prod_i H(X_i).
\end{array}
$$
The element $\displaystyle\prod_i \alpha_i \in \displaystyle\prod_i
H(X_i)$ is in the kernel of (1-shift), and hence there is an $\alpha \in
H(X)$ mapping to it.  By Yoneda, $\alpha$ corresponds to a natural
transformation
$$
\mbox{Hom}(-,X) \rarrow H,
$$
and the fact that $\alpha$ maps to $\prod \alpha_i \in H(\coprod X_i)$
means that the diagram

\begin{picture}(400,100)
\put(120,10) {$\mbox{Hom}(-,X)$}
\put(200,10) {$\longrightarrow$}
\put(260,10) {$H$}
\put(180,70) {$\mbox{Hom}(-,X_i)$}
\put(170,40) {$\swarrow$}
\put(240,40) {$\searrow$}
\end{picture}

\noindent commutes for all $i$.
\bigskip

\noindent{\bf Proof of (2).} It remains to show that
$$
\phi: \mbox{ Hom}(-,X) \rarrow H
$$
constructed above is an isomorphism.  Let us begin with objects $t \in T$.
We will show that, for all $t \in T$,
$$
\phi(t): \mbox{ Hom}(t,X) \rarrow H(t)
$$
is an isomorphism.

Observe the commutative diagram
$$
\begin{array}{ccccc}
&& \mbox{Hom}(t,X_0) && \\
& \swarrow && \searrow & \\
\mbox{Hom}(t,X) && \longrightarrow && H(t)
\end{array}
$$
Since we know that Hom$(t,X_0) \rarrow H(t)$ is surjective, it follows that
$$
\mbox{Hom}(t,X) \rarrow H(t)
$$
is surjective.  It remains only to prove it injective.

Let $f \in \mbox{ Hom}(t,X)$ with $\phi(t)(f) = 0$.  Now $f \in \mbox{ Hom}
(t,X) = \mbox{ Hom}\left(t, \begin{array}{c}
\mbox{hocolim} \\ \vspace{-.75cm} \\ \longrightarrow \end{array}
X_i\right)$.  But as $t \in T$ is compact, we have by
Lemma~\ref{maps to colimits} that
$$
\mbox{Hom}\left(t,\begin{array}{c} \mbox{hocolim} \\
\vspace{-.75cm} \\ \longrightarrow \end{array} X_i\right) =
\displaystyle\lim_\rightarrow \mbox{ Hom}
(t,X_i).
$$
In other words, there exists $f_i: t \rarrow X_i$ so that the composite
$$
t \stackrel{f_i}{\longrightarrow} X_i \longrightarrow X
$$
is $f$.  But the diagram

\begin{picture}(400,100)
\put(120,10) {$\mbox{Hom}(t,X)$}
\put(200,10) {$\longrightarrow$}
\put(205,1)  {$k$}
\put(260,10) {$H(t)$}
\put(180,70) {$\mbox{Hom}(t,X_i)$}
\put(170,40) {$j \swarrow$}
\put(240,40) {$\searrow$}
\end{picture}

\noindent
commutes, and $j(f_i) = f$ while $k(f) = 0$.  It follows that
$f_i \in \ker \{\mbox{Hom}(t,X_i) \rarrow H(t)\}$, that is, $(f_i,t) \in
U_{i+1}$.  This means that $f_i: t \rarrow X_i$ factors through the
map $h$ in the triangle

\begin{picture}(400,100)
\put(90,15) {$\displaystyle\left\{\coprod_{(f_{i},t) \in U_{i+1}}t\right\} =
K_{i+1}$}
\put(220,15) {$\longleftarrow$}
\put(225,4)  {$\scriptstyle{(1)}$}
\put(260,15) {$X_{i+1}$}
\put(210,70) {$X_i$}
\put(170,44) {$h \nearrow$}
\put(240,44) {$\searrow g$}
\end{picture}

\noindent
and hence $g \circ f_i = 0$.  But the map
$$
X_i \stackrel{g}{\longrightarrow} X_{i+1} \stackrel{\bar{g}}{\longrightarrow}
X
$$
satisfies
\begin{eqnarray*}
f & = & \{ \bar{g} \circ g\} \circ f_i \\
& = & \bar{g} \circ \{g \circ f_i\} \\
& = & 0.
\end{eqnarray*}
Thus $\phi(t): \mbox{ Hom}(t,X) \rarrow H(t)$ is an isomorphism whenever
$t \in T$.

Let ${\cal S} \subset {\cal T}$ be the full subcategory of objects $y \in {\cal
T}$
such that, for all $n\in \boldZ$, the map
$\phi(\Sigma^n y): \mbox{ Hom}(\Sigma^n y,X) \rarrow
H(\Sigma^n y)$ is an isomorphism.  Then the category ${\cal S}$ contains $T$,
and is closed with respect to the formation of coproducts and triangles.
To finish our proof of Theorem~\ref{Brown representability}, we need the lemma

\lem{$T$ generates} Let ${\cal S} \subset {\cal T}$
be a full, triangulated subcategory containing $T$ and closed under the
formation of ${\cal T}$-coproducts of its objects.  Then ${\cal S} =
{\cal T}$.\elem

\noindent{\bf Proof.}  Let ${\cal S}$ be the smallest subcategory of ${\cal T}$
which is full, triangulated, closed with respect to ${\cal T}$-coproducts
of its objects, and contains $T$.  Let $Z$ be an object of ${\cal T}$.
Let $H = \mbox{Hom}_{\cal T}(-,Z)$ be viewed as a homological functor
on ${\cal S}$.

Then ${\cal S}$ is compactly generated, with a generating set $T$.  We can
therefore apply what we have proved so far about Brown representability
to the functor $H$ on ${\cal S}$; there is an object $X$ of ${\cal S}$,
a natural transformation of functors on ${\cal S}$
$$
\phi: \mbox{Hom}_{\cal S} (-,X) \rarrow \mbox{Hom}_{\cal T}(-,Z),
$$
and this natural transformation is an isomorphism on a full, triangulated
subcategory of ${\cal S}$ containing $T$ and closed with respect to
${\cal S}$-coproducts of its objects.  But ${\cal S}$ is minimal with this
property; hence
$$
\mbox{Hom}_{\cal S} (-,X) \rarrow \mbox{Hom}_{\cal T} (-,Z)
$$
is an isomorphism of functors on ${\cal S}$.

By Yoneda's lemma, this means there is a morphism in ${\cal T}$
$$
X \rarrow Z
$$
so that, for every object, $s$ of ${\cal S}$,
$$
\mbox{Hom}(s,X) \rarrow \mbox{Hom}(s,Z)
$$
is an isomorphism.

Complete $X \rarrow Z$ to a triangle in ${\cal T}$
$$
\begin{array}{rcccc}
X && \longrightarrow && Z \\
& {\scriptstyle{(1)}} \nwarrow & & \swarrow & \\
&& Y &&
\end{array}
$$
One easily deduces that, for every object $s$ of $\cal S$, $\mbox{Hom}(s,Y)
= 0$.  But $T \subset {\cal S}$, and hence for every object $t \in T, \
\mbox{Hom}(t,Y) = 0$.  But because $T$ generates, $Y = 0$, and hence
$X \rarrow Z$ is an isomorphism.  Thus $Z$ is an object of ${\cal S}$,
and since $Z \in {\cal T}$ was arbitrary, ${\cal S} = {\cal T}$.
\hfill $\Box$

\section{The adjoint functor theorem and examples}
\label{S3}

\thm{adjoints} Let ${\cal S}$ be a compactly generated
triangulated category, ${\cal T}$ any triangulated category.
Let $F: {\cal S} \rarrow {\cal T}$ be a triangulated functor.  Suppose $F$
respects coproducts; the natural maps
$$
F(s_\lambda) \rarrow F\left(\displaystyle\coprod_{\lambda \in \Lambda}
s_\lambda\right)
$$
make $F\left(\displaystyle\coprod_{\lambda \in \Lambda}
s_\lambda\right)$ a coproduct of ${\cal T}$.  Then $F$ has a right
adjoint $G: {\cal T} \rarrow {\cal S}$.\ethm

\noindent{\bf Proof.}  Let $t$ be an object of ${\cal T}$, and
consider the functor on ${\cal S}$
$$
s \mapsto \mbox{ Hom}_{\cal T} (F(s),t).
$$
This functor is homological and takes coproducts to products; we
have
\begin{eqnarray*}
\mbox{Hom}_{\cal T}\left(F\left(\displaystyle\coprod_{\lambda \in \Lambda}
s_\lambda\right),t\right)
& = & \mbox{Hom}_{\cal T}\left(\coprod_{\lambda \in \Lambda}
F(s_\lambda),t\right) \\
& = & \prod_{\lambda \in \Lambda} \mbox{Hom}_{\cal T}(F(s_\lambda),t).
\end{eqnarray*}
Hence, by Theorem~\ref{Brown representability}, this functor is representable;
there is a $G(t) \in {\cal S}$ with
$$
\mbox{Hom}_{\cal T}(F(s),t) = \mbox{Hom}_{\cal S}(s,G(t)),
$$
and by standard arguments, $G$ extends to a functor, right adjoint to $F$.
\hfill $\Box$

\exm{$qc$ adjoints}  Let $f: X \rarrow Y$ be a
separated morphism of quasi-compact separated schemes. Then
$$
 Rf_*: D(qc/X) \rarrow D(qc/Y)
$$
has a right adjoint
$$
f^!:D(qc/Y) \rarrow D(qc/X).
$$\eexm

\noindent{\bf Proof.} We need only show that $Rf_*$ is triangulated and
respects coproducts; the fact that it is triangulated is obvious, the fact that
it respects coproducts is Lemma~\ref{$Rf_*$ respects sums}. \hfill $\Box$

\exm{$D$-module adjoints}  Let $f: X \rarrow Y$ be
a separated morphism of smooth, quasi-compact, separated schemes of finite
type over a field $k$ of characteristic 0.  Then
$$
Rf_+: D\left(\frac{qc \: D\mbox{-modules}}{X}\right) \rarrow D\left(\frac{qc \:
D\mbox{-modules}}{Y}\right)
$$
has a left adjoint.\eexm

\noindent{\bf Proof.}  Since $Rf_+$ is clearly a triangulated functor,
the real point is to prove that it respects coproducts.  But $Rf_+$
is given as
$$
Rf_+(x) \stackrel{\rm def}{=} Rf_*\left({\cal D}^{}_{Y \leftarrow X}
\otimes_{{\cal D}^{}_{X}}x\right)
$$
and tensor product trivially respects coproducts, while $Rf_*$ does by
Lemma~\ref{$Rf_*$ respects sums}. \hfill $\Box$
\bigskip

\exm{holomorphic adjoints}  With $f: X \rarrow
Y$ as in Example~\ref{$D$-module adjoints}, let $D(\mbox{holo}/X)$ and
$D(\mbox{holo}/Y)$ be as in Example~\ref{holo}; that is, $D(\mbox{holo}/X)$
is the smallest full, triangulated category of
$D\left(\frac{qc \: D\mbox{-modules}}{X}\right)$
closed with respect to coproducts and containing the bounded complexes
with holonomic cohomology.  It is well known that
$$
Rf_+: D\left(\frac{qc \: D\mbox{-modules}}{X}\right) \rarrow D\left(\frac{qc \:
D\mbox{-modules}}{Y}\right)
$$
takes complexes with holonomic cohomology to complexes with holonomic
cohomology.  Since $Rf_+$ is triangulated and respects coproducts, it
takes $D(\mbox{holo}/X)$ to $D(\mbox{holo}/Y)$.
It induces a functor, which we also denote $Rf_+$,
$$
Rf_+: D(\mbox{holo}/X) \rarrow D(\mbox{holo}/Y).
$$
This functor has a left adjoint.\eexm

\noindent{\bf Proof.}  By Example~\ref{holo}, $D(\mbox{holo}/X)$
is a compactly generated triangulated category.  The functor $Rf_+$ is the
restriction to $D(\mbox{holo}/X)$ of a functor on
$D\left(\frac{qc \: D\mbox{-modules}}{X}\right)$
respecting coproducts; hence it respects coproducts.  By Theorem~\ref{Brown
representability}, the adjoint exists.  \hfill $\Box$

\section{Commuting with coproducts}
\label{S4}

In Section~\ref{S3}, we constructed an adjoint to a functor $F: {\cal S}
\rarrow {\cal T}$.  Precisely, if $F: {\cal S} \rarrow {\cal T}$ is a
triangulated functor, one sometimes has a right adjoint $G: {\cal T}
\rarrow {\cal S}$.  Being a right adjoint, $G$ certainly respects products.
It turns out to be interesting to know when $G$ respects coproducts.

\thm{respecting coproducts} Let ${\cal S}$ be a compactly generated
triangulated category, and let ${\cal T}$ be any triangulated category.
Let $F: {\cal S} \rarrow {\cal T}$ be a triangulated functor respecting
coproducts, and let $G: {\cal T} \rarrow {\cal S}$ be its right adjoint,
which exists by Theorem~\ref{adjoints}.  Let $S$ be a generating set for
${\cal S}$.  Then $G: {\cal T} \rarrow {\cal S}$ respects coproducts if
and only if for every $s \in S, \ F(s)$ is a compact object of ${\cal T}$.\ethm

\noindent{\bf Proof.} $\Rightarrow$ Suppose $G$ preserves coproducts.  Let
$s \in S$ be some object.  Then

$$
\begin{array}{cccl}
\mbox{Hom}_{\cal T} \left(F(s),\displaystyle\coprod_{\lambda \in \Lambda}
x_\lambda\right) & = & \mbox{Hom}_{\cal S}\left(s,G\left(
\displaystyle\coprod_{\lambda \in \Lambda} x_\lambda\right)\right)& \\*[3pt]
& = & \mbox{Hom}_{\cal S}\left(s,\displaystyle\coprod_{\lambda \in \Lambda}
G\left(x_\lambda\right)\right)&\quad\hbox{because G commutes}  \\*[-10pt]
&&&\quad\hbox{with coproducts}\\*[10pt]
& = & \displaystyle\coprod_{\lambda \in \Lambda}\mbox{Hom}_{\cal S}\left(s,
G\left(x_\lambda\right)\right)&\quad\hbox{because $s$ is compact} \\*[20pt]
& = & \displaystyle\coprod_{\lambda \in \Lambda} \mbox{Hom}_{\cal T}
(F(s),x_\lambda),&
\end{array}
$$
and this proves that $F(s)$ is compact in ${\cal T}$.
\medskip

\nin
$\Leftarrow$ Suppose that, for all $s \in S, \ F(s)$ is compact.  Let
$\displaystyle\coprod_{\lambda \in \Lambda} x_\lambda$ be a coproduct
in ${\cal T}$.  Then for any $s \in S$,
$$
\begin{array}{cccc}
\mbox{Hom}\left(s,G\left(\displaystyle\coprod_{\lambda \in \Lambda}
x_\lambda\right)\right)
& = & \mbox{Hom}\left(F(s),\displaystyle\coprod_{\lambda \in \Lambda}
x_\lambda\right) &\\
& = & \displaystyle\coprod_{\lambda \in \Lambda} \mbox{Hom}
(F(s),x_\lambda)&\quad\hbox{because $F(s)$ is compact} \\*[20pt]
& = & \displaystyle\coprod_{\lambda \in \Lambda} \mbox{Hom}(s,G(x_\lambda))& \\
& = & \mbox{Hom}\left(s,\displaystyle\coprod_{\lambda \in \Lambda}
G(x_\lambda)\right)&\quad\hbox{because $s$ is compact}.
\end{array}
$$
In other words, the natural map
$$
\displaystyle\coprod_{\lambda \in \Lambda} G(x_\lambda) \rarrow
G\left(\displaystyle\coprod_{\lambda \in \Lambda} x_\lambda\right)
$$
induces a natural transformation
$$
\phi: \mbox{Hom}_{\cal S}
\left(-\:, \displaystyle\coprod_{\lambda \in \Lambda} G(x_\lambda)\right)
\rarrow \mbox{Hom}_{\cal S}\left(-\:,G\left(
\displaystyle\coprod_{\lambda \in \Lambda} x_\lambda\right)\right),
$$
and $\phi(s)$ is an isomorphism for all $s \in S$. But then in the triangle
$$
\displaystyle\coprod_{\lambda \in \Lambda} G(x_\lambda) \rarrow
G\left(\displaystyle\coprod_{\lambda \in \Lambda} x_\lambda\right)\la
Z\la \Sigma\left\{\displaystyle\coprod_{\lambda \in \Lambda}
G(x_\lambda)\right\}
$$
the object $Z$ must satisfy $Hom({\cal S}, Z)=0$. But as $\cal S$ generates,
$Z=0$ and
$$
\displaystyle\coprod_{\lambda \in \Lambda} G(x_\lambda) \rarrow
G\left(\displaystyle\coprod_{\lambda \in \Lambda} x_\lambda\right)
$$
is an isomoprphism.\hfill $\Box$
\bigskip

\exm{example}  Let $f: X \rarrow Y$ be a pseudo--coherent
proper morphism of
separated, quasi--compact schemes.  Suppose that $f$ is of finite
Tor-dimension.  Then $Rf_*: D(qc/X) \rarrow D(qc/Y)$ takes perfect complexes
to perfect complexes, by \cite{Ki}.

It follows that it takes a set of generators of $D(qc/X)$ to a set of
compact objects of $D(qc/Y)$.  Hence $f^!$ commutes with coproducts. What makes
this interesting is Theorem~\ref{$f^!$ commutes with sums}
which follows. But in the theorem we appeal to the projection formula,
and since I do not know a reference which proves it
in the generality we want, the following is the sketch
of a proof.\eexm

\pro{projection formula} Let $f:X\la Y$ be a morphism of separated,
quasi--compact schemes. Let $D(X)$ be the derived category of
all ${\cal O}_X$--modules, $D(Y)$ the derived category of all
${\cal O}_Y$--modules. Let $y$ be an object of $D(Y)$, $x$
an object of $D(X)$. Then there is a natural map, in $D(Y)$,
$$
y{^L\otimes}Rf_*x \la Rf_*\left\{Lf^*y{^L\otimes} x\right\}.
$$
If $y$ is in $D(qc/Y)\subset D(Y)$ and $x$ is in $D(qc/X)\subset
D(X)$, this natural map is an isomorphism.\epro

\noindent
{\bf Proof.}\ \ The existence and naturality of the map really comes from
the definition of $Rf_*$, $Lf^*$ and the derived tensor product.
To define $Rf_*$ of an object in $D(X)$, one expresses
the derived category $D(X)$ as
a quotient of the the homotopy category $K(X)$ by the acyclic
complexes $E(X)$, and notes that the subcategory $L(X)$ of Bousfield
local objects in $K(X)$ with respect to $E(X)$ maps isomorphically
to $D(X)$. Here, $L(X)$ consists of the special complexes of injectives
of Spaltenstein's \cite{S}. $Rf_*$ is just $f_*$ on $L(X)$.

The tensor product, and $Lf^*$, depend for their construction
on the fact that $D(X)$ is also isomorphic to the subcategory
of Bousfield colocal objects, denoted here $\tilde L(X)$. Concretely,
these are complexes of objects $j_!{\cal O}_U$, where $j:U\la X$
is the inclusion of an open affine, and $j_!$ is extension by zero.

Replacing $y$ by a colocal object on $Y$ and $x$ by a local
object on $X$, the tensor products become the natural ones
and we have an isomorphism
$$
y\otimes f_*x \simeq f_*\left\{f^*y\otimes x\right\}.
$$
The left hand side identifies immediately with $y{^L\otimes}Rf_*x$,
by definition of the derived tensor product and $Rf_*$. On
the right, the part in bracket identifies with
$Lf^*y{^L\otimes} x$, again by definition. Hence a natural isomorphism
$$
y{^L\otimes}Rf_*x \simeq f_*\left\{Lf^*y{^L\otimes} x\right\}.
$$
The problem is that, in general, $Lf^*y{^L\otimes} x$ is not
Bousfield local (=a complex of injectives), and hence $f_*$ of it
is not the same as $Rf_*$. But for any complex $x$, there is a natural
map $f_*x\la Rf_*x$, and this gives the natural map
$$
y{^L\otimes}Rf_*x \la Rf_*\left\{Lf^*y{^L\otimes} x\right\}.
$$
It remains to show that the restriction of this map to the
subcategories of complexes of quasi--coherents is an isomorphism.

Fix $x\in D(qc/X)$; we have a natural transformation of functors
in $y\in D(qc/Y)$.
First, the problem is now local in $Y$ and we may therefore assume
$Y$ affine. Secondly, on the category $D(qc/X)$, $Rf_*$ respects
coproducts by Lemma~\ref{$Rf_*$ respects sums}. Tensor product
and $Lf^*$ always respect sums, so the map is a natural transformation
of two functors in $y$ respecting coproducts. For each $y\in D(qc/Y)$,
let $\phi(y)$ be the map
$$
y{^L\otimes}Rf_*x \la Rf_*\left\{Lf^*y{^L\otimes} x\right\}.
$$
Let ${\cal R}\subset D(qc/Y)$ be the full subcategory of all $y$'s
such that $\phi(\Sigma^n y)$ is an isomorphism for all $n$.
The category $\cal R$ is closed with respect to triangles
and coproducts. It clearly contains ${\cal O}_Y$. Since $Y$
is affine, ${\cal O}_Y$ is ample and generates $D(qc/Y)$
by Example~\ref{quasi-projective variety}. By Lemma~\ref{$T$ generates},
it follows that $\cal R$ is all of $D(qc/Y)$.\qqed

\thm{$f^!$ commutes with sums} Let $f: X \rarrow Y$
be a morphism of schemes.  Suppose $Rf_*$ has a right adjoint $f^!$ which
commutes with coproducts.  Suppose $Y$ is quasi-compact and separated.
Then there is a natural isomorphism
of functors $D(qc/Y) \rarrow D(qc/X)$, which on objects gives
$$
f^!(y) \simeq (Lf^*y) \otimes_{{\cal O}_{X}} (f^!{\cal O}_Y).
$$
Conversely, if $f^!$ is naturally isomorphic to the functor on the right,
it respects coproducts.\ethm

\noindent{\bf Proof.} $\Leftarrow$ Suppose we have a natural isomorphism of
functors in $y$
$$
f^!(y) \simeq (Lf^*y) \otimes_{{\cal O}_{Y}}(f^!{\cal O}_Y).
$$
Since $Lf^*$ has a right adjoint, it respects coproducts.  Tensor
products respect coproducts, so we deduce that $f^!$ respects coproducts.
\medskip

\nin
$\Rightarrow$ We will show that
there is a natural map
$$
 (Lf^*y) \otimes_{{\cal O}_{X}}(f^!{\cal O}_Y)\la f^!(y)
$$
and that this map is an isomorphism whenever $y$ is compact. Then,
if $f^!$ respects coproducts, it will easily follow that
this natural map is an isomorphism for all $y$.

Let us prove a slightly more general fact. We will actually show that,
for any $y'$ in $D(qc/Y)$, there is a natural map
$$
 (Lf^*y) \otimes_{{\cal O}_{X}}(f^!y')\la f^!(y\otimes_{{\cal O}_Y}y')
$$
which is an isomorphism if $y$ is compact. The case $y'={\cal O}_Y$
is then the above.

In any case, there is a natural map
$$
\mu: Rf_* f^! y' \rarrow y',
$$
the counit of adjunction.  For every $y \in D(qc/Y)$,
$$
Rf_*\left[(Lf^*y) \otimes_{{\cal O}_{X}} f^!y'\right]
= y \otimes_{{\cal O}_{Y}} Rf_* f^! y'
$$
by the projection formula. Hence there is a natural map
$$
\begin{array}{rcc}
Rf_*\left[Lf^*(y) \otimes_{{\cal O}_{X}} f^!y'\right]
& = & y \otimes_{{\cal O}_{Y}} Rf_*f^!y' \\*[10pt]
&& \hspace*{1cm} \downarrow 1_y \otimes \mu \\*[10pt]
&& y \otimes_{{\cal O}_{Y}} y'.
\end{array}
$$
By adjunction, we have a map
$$
Lf^*(y) \otimes_{{\cal O}_{X}} f^!y'\la f^!\left(y\otimes_{{\cal
O}_{Y}}y'\right)
$$
This is our natural map. We need to show that for compact $y$
it is an isomorphism.

Let $y$ be compact, and let $x$ be an arbitrary complex in $D(qc/X)$.
It suffices to show that the natural map above induces an isomorphism
after applying $Hom(x,-)$. Let us therefore reflect what $Hom(x,-)$ does.
To begin with, put $\hat y={\cal RH}{\rm om}_{{\cal O}_Y}(y,{\cal O}_Y)$.
Then $Lf^*\hat y={\cal RH}{\rm om}_{{\cal O}_X}(Lf^*y,{\cal O}_X)$, and since
$y$ and $Lf^*y$ are perfect complexes, there are natural isomorphisms
$$
Hom_X(-\otimes Lf^*\hat y\,,\,-)=Hom_X(-\,,\,Lf^*y\otimes-)
$$
and
$$
Hom_Y(-\otimes\hat y\,,\,-)=Hom_Y(-\,,\,y\otimes-).
$$
Now, the map
$$
Lf^*(y) \otimes_{{\cal O}_{X}} f^!y'\la f^!\left(y\otimes_{{\cal
O}_{Y}}y'\right)
$$
induces a map after applying $Hom(x,-)$. By definition, this takes
a map $\gamma\in Hom\left(x\,,\,Lf^*(y) \otimes_{{\cal O}_{X}} f^!y'\right)$ to
a map $x\la f^!(y\otimes_{{\cal O}_{Y}} y')$. By the adjunction, this
map corresponds to a map $\gamma':Rf_*x\la y\otimes_{{\cal O}_{Y}} y'$.
But we  know what $\gamma'$ is; it is the composite
of $Rf_*\gamma$ with the natural counit
$$
\begin{array}{rcc}
Rf_*\left[Lf^*(y) \otimes_{{\cal O}_{X}} f^!y'\right]
& = & y \otimes_{{\cal O}_{Y}} Rf_*f^!y' \\*[10pt]
&& \hspace*{1cm} \downarrow 1_y \otimes \mu \\*[10pt]
&& y \otimes_{{\cal O}_{Y}} y'.
\end{array}
$$
We need to show this correspondence to be an isomorphism. But
\begin{eqnarray*}
\mbox{Hom}_{{\cal O}_{X}}
\left(x\,,\,Lf^*y \otimes_{{\cal O}_{X}} f^!y'\right)
& = & \mbox{Hom}_{{\cal O}_{X}}\left(x \otimes Lf^*\hat y\,,\,
f^!y'\right) \\
& = & \mbox{Hom}_{{\cal O}_{Y}}\left(Rf_*
\{x \otimes Lf^*\hat y\}\,,\,y'\right) \\
& = & \mbox{Hom}_{{\cal O}_{Y}}\left(
Rf_*(x)\otimes \hat y\,,\,y'\right) \\
& = & \mbox{Hom}_{{\cal O}_Y}
\left(Rf_*(x)\,,\,y \otimes_{{\cal O}_{Y}} y'\right),
\end{eqnarray*}
where the third equality is the projection formula.
This isomorphism is easily identified with $Hom(x,-)$ applied to
the map
$$
\phi:Lf^*(y) \otimes_{{\cal O}_{X}} f^!y'\la f^!\left(y\otimes_{{\cal
O}_{Y}}y'\right).
$$
We therefore know that $\phi$ is an isomorphism if $y$ is compact.

Now assume that $f^!$ respects coproducts. For fixed $y'$, the functor
$$
y \mapsto (Lf^*y) \otimes_{{\cal O}_{X}} f^!y'
$$
is a triangulated functor in $y$; $Lf^*$ is, as is tensor product.  The
functor $f^!$ is the adjoint of a triangulated functor, hence triangulated;
see \cite{SH}, Lemma 3.9.  Thus $\phi$ is a natural transformation of
triangulated
functors, both of which respect coproducts.
Let ${\cal S} \subset D(qc/Y)$ be the full subcategory
$$
{\cal S} = \left\{y \in {\cal O}b\left[D(qc/Y)\right]| \: \phi(\Sigma^ny)
\mbox{ is
an isomorphism for all } n \in \boldZ\right\}.
$$
Then ${\cal S}$ contains the generating set of compact objects,
is triangulated and closed
with respect to $D(qc/Y)$-coproducts.  Thus ${\cal S} = D(qc/Y)$.
Hence, $\phi$ is a natural isomorphism. The theorem is the special
case $y'={\cal O}_Y$ of the above.  \hfill $\Box$
\bigskip

\rmk{traditional} In the traditional literature
on the subject, $f^!{\cal O}_Y$ is called the {\em dualizing complex},
and plays a key role in the theory.\ermk

\section{A sheaf version}
\label{S5}

The traditional way to state Grothendieck's duality theorem comes in a sheaf
version.  Let $f: X \rarrow Y$ be a proper morphism of noetherian, separated
schemes.  One would like to deduce that $f^!$, the adjoint we have for
$Rf_*$, gives an isomorphism in the category of sheaf homomorphisms
$$
{\cal RH}om(Rf_*x,y) \simeq Rf_* {\cal RH}om(x,f^!y).
$$
Note that the counit $u: Rf_*f^!y \rarrow y$ defines in any case a well-defined
natural transformation
$$
\phi: Rf_* {\cal RH}om(x,f^!y) \rarrow  {\cal RH}om(Rf_*x,y),
$$
and the only question is whether it is an isomorphism.  Once we apply the
functor $H^0(Y,-)$, it becomes an isomorphism; this is because $f^!$ is adjoint
to $Rf_*$.  To say that the map is an isomorphism of sheaves is to say that
the derived functor of $\Gamma(U,-)$ gives an isomorphism $R\Gamma(U,\phi)$
for every open set $U \subset Y$.  Concretely, it says that if we take
the commutative diagram
$$
\begin{array}{ccc}
f^{-1}U & \stackrel{j'}{\hookrightarrow} & X \\
f' \downarrow && \downarrow f \\
U & \stackrel{j}{\hookrightarrow} & Y,
\end{array}
$$
then $(f')^!j^*$ and $(j')^*f^!$ are naturally isomorphic.

It is easy to see that $j^*Rf_* = Rf'_*(j')^*$.  Taking right adjoints,
we deduce

\lem{$^*$ OK} If the diagram
$$
\begin{array}{ccc}
f^{-1}U & \stackrel{j'}{\hookrightarrow} & X \\
f' \downarrow && \downarrow f \\
U & \stackrel{j}{\hookrightarrow} & Y,
\end{array}
$$
is given by pulling back an open immersion $j:U\la Y$, then
there is a natural isomorphism
$$f^!R{j_*} = R{j'_*}(f')^!$$.\qqed\elem

Therefore
$$
(j')^*\{Rj'_*(f')^!\}j^* = (j')^*\{f^!R{j_*}\}j^*.
$$
Hence $(f')^!j^* = (j')^*f^!R{j_*}j^*$ because $(j')^*R{j'_*}$ is the
identity functor.  Thus we must convince ourselves that the natural
unit of adjunction $1 \rarrow R{j_*}j^*$ induces an isomorphism
$$
(j')^*f^! \rarrow (j')^*f^!R{j_*}j^*.
$$
Let $Z=Y-U$ be the (closed) complement of
$U\subset Y$.  There is a triangle of functors
$$
\begin{array}{c}
R\Gamma_Z \rarrow 1 \rarrow Rj_*j^* \\ \longleftarrow \\ (1)
\end{array}
$$
where $R\Gamma_Z$ is Grothendieck's local cohomology functor. Sometimes
$R\Gamma_Z$ is denoted $i_*i^!$, because it is also the counit
of an adjunction. But since the $i^!$ and $i_*$ are not the type
of map we have been considering here, the notation might lead to
confusion. Note that $Z\subset Y$ is a Zariski closed subset, but we
have {\em not} given it a scheme structure.  Hence the $i_*$ and $i^!$
of this article make no sense.

It suffices to show that
$$
(j')^* f^! R\Gamma_Z
$$
is the zero functor.  It is enough to show this for open sets $U \subset X$
which form a basis for the topology.  This is what we will do.
\bigskip

\pro{the vanishing} Let $f: X \rarrow Y$ be a
morphism of schemes such that $Rf_*: D(qc/X) \rarrow D(qc/Y)$ has a right
adjoint $f^!$.  Suppose $f^!$ respects coproducts.  Suppose $U \subset Y$
is an open subset, $Z = Y-U$ the complement. Suppose $U$ and $Y$
are quasi--compact and separated.
Then, in the notation above, the composite
$$
(j')^* f^! R\Gamma_Z = 0.
$$\epro

\noindent{\bf Proof.}\ \
We wish to show that $(j')^* f^!$ vanishes on any
object of the form $R\Gamma_Z (y)$. This means concretely
that if $y$ is a complex which is acyclic away from
$Z$, then $f^!y$ must be shown acyclic off $f^{-1}Z$.
But by
Theorem~\ref{$f^!$ commutes with sums},
$$
f^!y=Lf^*y\otimes f^!{\cal O}_Y.
$$
Clearly, if $y$ is supported on $Z$, then $Lf^*y$ is
supported on $f^{-1}Z$, and hence so is
its tensor product with $f^!{\cal O}_Y$.\qqed

 From the work of Verdier for the Noetherian case, Lipman in
general, we know that for the bounded--below derived category
more is true. Let us state their theorem, then prove it by
coproduct techniques.

\pro{the vanishing, bounded} Let $f: X \rarrow Y$ be a
pseudo--coherent,
proper morphism of quasi--compact, separated
schemes. Then $Rf_*: D^+(qc/X) \rarrow D^+(qc/Y)$ has a right
adjoint $f^!$.
Furthermore, in the notation above, the composite
$$
(j')^* f^! R\Gamma_Z = 0.
$$\epro

\noindent{\bf Proof.}\ \
Let us begin by showing that the right
adjoint of $f^!: D(qc/Y) \rarrow D(qc/X)$ takes
bounded--below complexes to bounded--below complexes, and
is therefore also a right adjoint to
$Rf_*: D^+(qc/X) \rarrow D^+(qc/Y)$. Since $X$ is quasi-compact,
it may be covered by finitely many open affines, say $n$ of them.
But since $f$ is separated, the open affines may be used to
compute $Rf_*$, via the \v Cech complex. Since the \v Cech
complex has only $(n+1)$ terms, it follows that if $x\in D^-(qc/X)$
vanishes above dimension $l$, then $Rf_*x$ vanishes above dimension
$l+n$. In the notation of {\it t}--structure truncations,
if $x\in {D(qc/X)}^{\le l}$, then $Rf_*x\in {D(qc/Y)}^{\le l+n}$.

Pick any $x\in {D(qc/X)}^{\le l}$, and $y\in {D(qc/Y)}^{\ge l+n+1}$.
Then
$$
Hom(x,f^!y)=Hom(Rf_*x,y)=0
$$
since $Rf_*x\in {D(qc/Y)}^{\le l+n}$ and  $y\in {D(qc/Y)}^{\ge l+n+1}$.
But $x\in {D(qc/X)}^{\le l}$ was arbitrary; thus
$f^!y\in {D(qc/X)}^{\ge l+1}$. In other words, $y\in {D(qc/Y)}^{\ge l}$
implies $f^!y\in {D(qc/X)}^{\ge l-n}$. In particular, if $y\in
D^+(qc/Y)$, then $f^!y\in D^+(qc/X)$.

It remains  to show the vanishing of
$
(j')^* f^! R\Gamma_Z.
$
We need to show that if $y$ is an object of $D^+(qc/Y)$,
which is acyclic off $Z$, then $f^!y$ is acyclic
off $f^{-1}Z$.
Let us first make an observation. Suppose
$y$ is arbitrary, vanishing off $Z$.
Suppose $Y=V_1\cup V_2$ expresses $Y$ as a union
of two open sets $V_1$ and $V_2$. There is then a triangle
$$
y\la {j^{}_{V_1}}_*{j^{}_{V_1}}^*y\oplus
{j^{}_{V_2}}_*{j^{}_{V_2}}^*y
\la {j^{}_{V_1\cap V_2}}_*{j^{}_{V_1\cap V_2}}^*y\la\Sigma y
$$
where $j^{}_{W}:W\hookrightarrow Y$ is the inclusion of the open set
$W$ in $Y$. From the triangle, it clearly suffices to show the
vanishing of $(j')^* f^!$ on the complexes
${j^{}_{W}}_*{j^{}_{W}}^*y$ where $W$ is any of $V_1$, $V_2$
or $V_1\cap V_2$. Suppose $X$ can be covered by $n$ affines.
If $V_1$ is affine and $V_2$ the union of $n-1$ affines, then
each of $V_1$, $V_2$ and $V_1\cap V_2$ can be covered by at most
$n-1$ affines. By induction we therefore easily show that it
suffices to prove the vanishing of $(j')^* f^!$ on
${j^{}_{W}}_*{j^{}_{W}}^*y$ where $W\subset X$ is affine;
in other words, we need to show that
$f^!{j^{}_{W}}_*{j^{}_{W}}^*y$ vanishes off $f^{-1}Z$.

Even better, we may replace $y$ by $z={j^{}_{W}}^*y$. It
will suffice to show that, given an open affine $W\subset Y$
and
 $z\in D^+(qc/W)$ whose support is in $Z\cap W$,
then $f^!{j^{}_{W}}_*z$ is supported on $f^{-1}Z$.
Now consider the pullback square
$$
\begin{array}{ccc}
f^{-1}W & \stackrel{j_{f^{-1}W}}{\hookrightarrow} & X \\
f' \downarrow && \downarrow f \\
W & \stackrel{j_W}{\hookrightarrow} & Y.
\end{array}
$$
By Lemma~\ref{$^*$ OK} there is a natural isomorphism
${\{j_{f^{-1}W}\}}_*{\{f'\}}^!z=f^!{\{j_W\}}_*z$. We
need to show that this complex (either of the two isomorphic
versions) is acyclic off $f^{-1}Z$. From the description as
${\{j_{f^{-1}W}\}}_*{\{f'\}}^!z$, it clearly suffices to
show that
${\{f'\}}^!z$ is supported on $f^{-1}\{Z\cap W\}$. In
other words, we are reduced to studying
the problem for the map $f:f^{-1}W\la W$. Thus we may assume $Y$
affine.

Next it is clear that if $Z=\cap Z_i$, it is
enough to prove the statement for each $Z_i$. We may therefore
assume that $Z$ is a ``divisor'' in $W$.
That is, there exists a global function $\gamma\in\Gamma(Y,{\cal O}_Y)$,
so that $Z$ is the divisor defined by $\gamma$.

On $Y$ we have perfect complexes, namely the
desuspensions of the mapping cones
of the maps
$$
\gamma^k:{\cal O}_Y\la {\cal O}_Y.
$$
Call these perfect complexes $b_k$. Then $b_k$ fits in a triangle
$$
b_k\la {\cal O}_Y\stackrel{\gamma^k}\la {\cal O}_Y \la \Sigma b_k.
$$
There is a map
$b_k\la b_{k+1}$, given by completing the following commutative
square
$$
\begin{array}{ccc}
{\cal O}_Y  &  \stackrel{\gamma^k}\la  &   {\cal O}_Y  \\
1\downarrow  &                       &  \downarrow\gamma  \\
{\cal O}_Y  &  \stackrel{\gamma^{k-1}}\la  &  {\cal O}_Y
\end{array}
$$
to a map of triangles
$$
\begin{array}{ccccccc}
b_k&\la     &{\cal O}_Y  &  \stackrel{\gamma^k}\la  &   {\cal O}_Y   &\la
&\Sigma b_k\\
\downarrow& &1\downarrow  &                       &  \downarrow\gamma &
&\downarrow \\
b_{k+1}&\la &{\cal O}_Y  &  \stackrel{\gamma^{k+1}}\la  &  {\cal O}_Y
&\la &\Sigma b_{k+1}
\end{array}
$$
There is furthermore a map $b_k\la {\cal O}_Y$, which is
a map from the direct system.

Applying the functor $(-)\otimes z$ to this direct system, we
get a direct system with a map
$$
b_k\otimes z\la {\cal O}_Y\otimes z=z.
$$
This induces a (non--canonical) map on homotopy colimits
$$
\hbox{hocolim}\left(b_k\otimes z\right)\la z
$$
and since the cohomology of $z$ is supported on $Z$, that
is annihilated by some power of $\gamma$, it is easy to
show that the map is a cohomology isomorphism, hence an
isomorphism in the derived category.

Note that in the construction above, $b_k\otimes z$
can be obtained by tensoring the triangle
$$
b_k\la {\cal O}_Y\stackrel{\gamma^k}\la {\cal O}_Y \la \Sigma b_k
$$
with the object $z$. There is a triangle
$$
b_k\otimes z\la z\stackrel{\gamma^k}\la z \la \Sigma b_k\otimes z.
$$
We know that $z\in D^+(qc/Y)$. Suppose $z\in D(qc/Y)^{\ge l}$.
Then from the triangle, $b_k\otimes z$ also lies in $D(qc/Y)^{\ge l}$.

Because $z$ is the homotopy colimit of $b_k\otimes z$, there is a triangle
on $Y$
$$
\bigoplus_{k}\left[ b_k\otimes z\right]\la
\bigoplus_{k}\left[ b_k\otimes z\right]\la
z\la
\Sigma\bigoplus_{k}\left[ b_k\otimes z\right]
$$
which expresses $z$ as the homotopy colimit.
Applying $f^!$, we have a triangle
$$
f^!\left\{\bigoplus_{k}\left[ b_k\otimes z\right]\right\}\la
f^!\left\{\bigoplus_{k}\left[ b_k\otimes z\right]\right\}\la
f^!z\la
f^!\left\{\Sigma\bigoplus_{k}\left[ b_k\otimes z\right]\right\}
$$
and to show that $f^!z$ is supported on $f^{-1}Z$,
it suffices to show that the other two terms in the triangle
are. Now note that in any case we there is a natural isomorphism
$$
f^!\left[ b_k\otimes z\right]=
Lf^* b_k\otimes f^!z
$$
because $b_k$ is compact, and by the proof
Theorem~\ref{$f^!$ commutes with sums}. Now
$Lf^* b_k$ is supported on $f^{-1}Z$,
hence its tensor product with $f^!{j^{}_W}_*z$
also is. It follows that
$f^!\left[ b_k\otimes z\right]$
is supported on $f^{-1}Z$. We need to show that
$\displaystyle
f^!\left\{\bigoplus_{k}\left[ b_k\otimes z\right]\right\}$
is supported on $f^{-1}Z$. It will clearly suffice to show that
the natural map
$$
\bigoplus_{k}f^!\left[ b_k\otimes z\right]\la
f^!\left\{\bigoplus_{k}\left[b_k\otimes z\right]\right\}
$$
is an isomorphism. This we will now do.

Now let $x$ be an arbitrary perfect complex on $X$. Since we
are only assuming that the map $f:X\la Y$ is proper and of finite type,
we do not know that $Rf_*x$ is perfect. However, we do know, by
\cite{Ki},
that locally it can be resolved by finite dimensional vector bundles
to arbitrary length. That is, $Y$ can be covered by open sets
$V$, and for each $V$ there is a triangle
$$
q\la p\la Rf_*x\la\Sigma q
$$
where $p$ is perfect and $q\in D(qc/V)^{\le l-1}$.
Since $Y$ is affine, this can even be done globally.
Such a triangle exists on all of $Y$.
We deduce
\begin{eqnarray*}
\hbox{Hom}\left(x,f^!\left\{\bigoplus_{k}\left[b_k\otimes z\right]
\right\}\right)
&=&\hbox{Hom}
\left(Rf_*x,\bigoplus_{k}\left[b_k\otimes z\right]\right) \\
&=&\hbox{Hom}
\left(p,\bigoplus_{k}\left[b_k\otimes z\right]\right) \\
&=&\bigoplus_{k}\hbox{Hom}
\left(p,\left[b_k\otimes z\right]\right) \\
&=&\bigoplus_{k}\hbox{Hom}
\left(Rf_*x,\left[b_k\otimes z\right]\right) \\
&=&\bigoplus_{k}\hbox{Hom}
\left(x,f^!\left[b_k\otimes z\right]\right)\\
&=&\hbox{Hom}
\left(x,\bigoplus_{k}f^!\left[b_k\otimes z\right]\right)
\end{eqnarray*}
But we have a map, defined on $Y$,
$$
\bigoplus_{k}f^!\left[b_k\otimes z\right]\la
f^!\left\{\bigoplus_{k}\left[b_k\otimes z\right]\right\}
$$
and we have just shown that if we apply the functor
$Hom(x,-)$ to this map where $x\in D(qc/X)$ is compact,
we get an isomorphism. But then the map
$$
\hbox{Hom}\left(x,
\bigoplus_{k}f^!\left[b_k\otimes z\right]
\right)\la
\hbox{Hom}\left(x,
f^!\left\{\bigoplus_{k}\left[b_k\otimes z\right]\right\}
\right)
$$
is
an isomorphism for all $x$ in the subcategory generated by the
compacts in $D(qc/X)$; that is for any $x\in D(qc/X)$. Thus
the map
$$
\bigoplus_{k}f^!\left[b_k\otimes z\right]\la
f^!\left\{\bigoplus_{k}\left[b_k\otimes z\right]\right\}
$$
is an isomorphism in $D(qc/X)$.\qqed

\rmk{Verdier} The results of Verdier apply not only
to open immersions, but to arbitrary flat map. In other words, he proves,
in Theorem~2 of \cite{V}, the following.
Suppose we have a cartesian square of
noetherian schemes

$$
\begin{array}{ccc}
X'&\stackrel {g'}\longrightarrow& X\\
{\scriptstyle f'}\!\downarrow\,\quad &   &{\scriptstyle
f}\!\!\downarrow\,\,\,\, \\
Y'&\stackrel {g}\longrightarrow& Y
\end{array}
$$

\nin
with $f$ and $f'$ proper, $g$ and $g'$ flat. Then there is a natural
isomorphism

$$
\{g'\}^*f^!=\{f'\}^!g^*.
$$
What interests us here is not so much the best statement possible,
but the relation with preserving coproducts. To illustrate this,
we will give the following counterexample.\ermk

\exm{counterexample} Let $R$ be a noetherian ring (eg $\Bbb Z$), and
let $S=R[\epsilon]/(\epsilon^2)$. There is a homomorphism
$S\la R$ sending $\epsilon$ to 0. This gives a map of schemes
from $X=Spec(R)$ to $Y=Spec(S)$. This map is certainly proper,
and of finite type. Let us denote this map $f:X\la Y$.

For affine maps, $f^!$ is easy to describe. For any $S$--module $N$
and $R$--module $M$, we can view $M$ as an $S$--module
via the homomorphism $S\la R$. This is the functor
$f_*:\{R\!-\!\hbox{mod}\}\la
\{S\!-\!\hbox{mod}\}$. There is a
canonical ismorphism
$$
Hom_S\left(M,N\right)=Hom_R\left(M,Hom_S(R,N)\right).
$$
This allows us to view $f^!$ as the derived $Hom$
$$
f^!N=RHom_S(R,N)
$$
for any $N\in D(qc/Y)=D(S)$. One way to get the derived functor
of $Hom$ is to take projective resolutions in the first variable.
There is an obvious projective resolution for $R$ as an $S$--module; the
chain complex
$$
\cdots\stackrel\epsilon\la S\stackrel\epsilon\la S\la R\la 0
$$
is exact. Let $N$ be $R$, viewed as an $S$--module. Then
$$
f^!R= RHom_S(R,R)=\prod_{i=0}^\infty \Sigma^{-i}R
$$
is the complex which is $R$ in every positive dimension, with
differential zero.

Now consider the complex
$$
N=\prod_{k=0}^\infty\Sigma^{k}R
$$
that is, the complex with zero differential which is $R$ in every negative
dimension, viewed as a complex of $S$--modules. Because $f^!$
is a right adjoint, it respects products. Thus,
\begin{eqnarray*}
f^!N&=&f^!\left\{\prod_{k=0}^\infty\Sigma^{k}R\right\}\\
&=& \prod_{k=0}^\infty f^!\left\{\Sigma^{k}R\right\}\\
&=& \prod_{k=0}^\infty\prod_{i=0}^\infty \Sigma^{k-i}R
\end{eqnarray*}
Choose an element $\gamma\in R\subset S=R[\epsilon]/(\epsilon^2)$.
If we restrict to the open subset where $\gamma$ is inverted,
then $f^!N$ restricts to
$$
{\{j'\}}^*f^!N=\left\{f^!N\right\}\otimes_R
R\left[\frac1\gamma\right].
$$
On the other hand, if we first restrict to the open subset,
note that we are in exactly the same situation; $S\left[\frac1\gamma\right]
=R\left[\frac1\gamma\right][\epsilon]/(\epsilon^2)$. The
complex $N$ restricts to
$$
j^*N=\prod_{k=0}^\infty \Sigma^{k}R\left[\frac1\gamma\right].
$$
If this is not clear, note that because the $R$'s are placed
in different dimensions, the coproduct agrees with the product.
The natural map
$$
\bigoplus_{k=0}^\infty \Sigma^{k}R\la
\prod_{k=0}^\infty \Sigma^{k}R=N
$$
is a homology isomorphism. But $j^*$, being
a left adjoint, respects coproducts; hence
$$
j^*\left(\bigoplus_{k=0}^\infty \Sigma^{k}R\right)=
\bigoplus_{k=0}^\infty j^*\Sigma^{k}R
=\bigoplus_{k=0}^\infty \Sigma^{k}R\left[\frac1\gamma\right]
$$
and, once again, the natural map
$$
\bigoplus_{k=0}^\infty \Sigma^{k}R\left[\frac1\gamma\right]\la
\prod_{k=0}^\infty \Sigma^{k}R\left[\frac1\gamma\right]
$$
is a homology isomorphism. We can therefore use our last
computation, replacing $R$ by $R\left[\frac1\gamma\right]$, to deduce
that
$$
{\{f'\}}^!j^*N=\prod_{k=0}^\infty\prod_{i=0}^\infty
\Sigma^{k-i}R\left[\frac1\gamma\right]
$$

In other words, we have just computed everything. It remains to check
whether the natural map is an isomorphism
${\{j'\}}^*f^!N\la {\{f'\}}^!j^*N$. The map is just
$$
 \left\{\prod_{k=0}^\infty\prod_{i=0}^\infty \Sigma^{k-i}R\right\}
\otimes_R
R\left[\frac1\gamma\right]\la
\prod_{k=0}^\infty\prod_{i=0}^\infty
\Sigma^{k-i}R\left[\frac1\gamma\right]
$$
If we just look at the induced map on $H^0$, the product is over
all $k=i$. The map on $H^0$ is therefore
$$
\left\{\prod_{i=0}^\infty R\right\}\otimes_R
R\left[\frac1\gamma\right]\la
\prod_{i=0}^\infty R\left[\frac1\gamma\right]
$$
This map clearly is not an isomorphism in general. If $\gamma$ is
neither nilpotent nor
invertible in $R$, then the element
$$
\prod_{i=0}^\infty \frac1{\gamma^i}
$$
is a well--defined member of the right hand side, but not the image of
something on the left.\eexm

\end{document}